\documentclass{amsart} 

\usepackage{verbatim}
\usepackage{amsmath}
\usepackage{amsfonts}
\usepackage{amssymb}
\usepackage{amscd}
\usepackage{latexsym}

\newcommand{\gld}{\operatorname{gldim}}
\newcommand{\pd}{\operatorname{pd}}

\newcommand{\dep}{\operatorname{depth}}

\newcommand{\reg}{\operatorname{reg}}

\newcommand{\lm}{\lim_{n \ra \infty}}
\newcommand{\cat}{\operatorname{-dist}}
\newcommand{\Cat}{\operatorname{-Dist}}

\newcommand{\sat}{\operatorname{sat}}

\newcommand{\spec} {\mbox{Spec}\ }
\newcommand{\ext} {\operatorname{Ext}}
\newcommand{\caphom} {\operatorname{Hom}}
\newcommand{\ann} {\operatorname{ann}}
\newcommand{\Llra}{\Longleftrightarrow}
\newcommand{\lra} {\longrightarrow}

\newcommand{\da} {\downarrow}

\newcommand{\ra}{\rightarrow}

\newcommand{\D}{\displaystyle}
\newcommand{\mc}{\mathcal}
\newcommand{\mf}{\mathfrak}
\newcommand{\cha}{\operatorname{char}}
\newcommand{\mb}{\mathbb}

\newcommand{\pgl}{\operatorname{PGL}}
\newcommand{\gl}{\operatorname{GL}}
\newcommand{\m}{\mf{m}}

\newcommand{\diag}{\operatorname{diag}}
\newcommand{\uext}{\underline{\operatorname{Ext}}}

\newcommand{\uhom}{\underline{\operatorname{Hom}}}
\renewcommand{\hom}{\operatorname{Hom}}

\newcommand{\GK}{\operatorname{GK}}

\newcommand{\wt}{\widetilde}
\newcommand{\ndo}{\operatorname{End}}
\newcommand{\coH}{\operatorname{H}}
\newcommand{\ucoH}{\underline{\operatorname{H}}}

\newcommand{\tor}{\operatorname{Tor}}
\newcommand{\utor}{\underline{\operatorname{Tor}}}

\newcommand{\im}{\operatorname{Im}}
\newcommand{\cd}{\operatorname{cd}}

\newcommand{\Proj}{\operatorname{-Proj}}
\newcommand{\proj}{\operatorname{-proj}}
\newcommand{\rprojnodash}{\operatorname{proj}}
\newcommand{\coh}{\operatorname{coh}}
\newcommand{\Gr}{\operatorname{-Gr}}

\newcommand{\gr}{\operatorname{-gr}}
\newcommand{\Qgr}{\operatorname{-Qgr}}
\newcommand{\qgr}{\operatorname{-qgr}}
\newcommand{\Tors}{\operatorname{-Tors}}
\newcommand{\tors}{\operatorname{-tors}}

\newcommand{\aut}{\operatorname{Aut}}

\newcommand{\regmult}{\circ} 
\newcommand{\twistmult}{\star} 
\newcommand{\Tmod}[1]{T^{#1}}
\newcommand{\pointmod}[1]{P({#1})}
\newcommand{\specmod}[1]{M^{#1}}

\theoremstyle{plain}
\newtheorem{theorem}{Theorem}[section]
\newtheorem{lemma}[theorem]{Lemma}

\newtheorem{corollary}[theorem]{Corollary}
\newtheorem{proposition}[theorem]{Proposition}

\theoremstyle{definition}
\newtheorem{definition}[theorem]{Definition}
\newtheorem{notation}[theorem]{Notation}
\newtheorem{example}[theorem]{Example}

\newtheorem{hypothesis}[theorem]{Standing Hypothesis}

\theoremstyle{remark}

\numberwithin{equation}{section}
\makeatletter                        
\let\c@equation\c@theorem              
\makeatother     

\begin{document}

\title{Generic Noncommutative Surfaces}
\author{Daniel Rogalski}
\address{University of Washington, Department of Mathematics, Box 354350, Seattle WA 98195, USA.} 
\email{rogalski@math.washington.edu}


\thanks{This research was supported in part by NSF
 grants NSF-DMS-9801148 and NSF-DMS-0202479, as well as 
a Clay Mathematics Institute Liftoff Fellowship. \hfill}

\abstract We study a class of noncommutative surfaces, and their higher dimensional 
analogues, which come from generic subalgebras of twisted homogeneous coordinate 
rings of projective space.  Such rings provide answers to several open questions in 
noncommutative projective geometry.  Specifically, these rings $R$ are the first 
known graded algebras over a field $k$ which are noetherian but not strongly 
noetherian: in other words, $R \otimes_k B$ is not noetherian for some choice of 
commutative noetherian extension ring $B$.  This answers a question of Artin, Small, 
and Zhang.  The rings $R$ are also maximal orders, but they do not satisfy all of the 
$\chi$ conditions of Artin and Zhang.  In particular, they satisfy the $\chi_1$ 
condition but not $\chi_i$ for $i \geq 2$, answering a question of Stafford and Zhang 
and a question of Stafford and Van den Bergh.  Finally, we show that the noncommutative
scheme $R\proj$ has finite global dimension.  
\endabstract

\maketitle
%\tableofcontents

\section{Introduction}

In algebraic geometry, the correspondence between projective schemes over a field $k$ 
and $\mb{N}$-graded $k$-algebras is classical.  In the last decade or so, many of the 
ideas of projective geometry have been successfully generalized to the setting of 
noncommutative graded rings.  This new subject is known as noncommutative 
projective geometry, and while of theoretical interest in its own right, it has also 
provided the solutions to many purely ring-theoretic questions.  For example, the 
graded domains of dimension two, which correspond to noncommutative curves, have now 
been completely classified \cite{ASt}.   The noncommutative analogs of the projective 
plane $\mb{P}^2$ have been identified and classified as well. The generic 
noncommutative projective plane is called the Sklyanin algebra; despite its 
simple presentation by generators and relations, before the geometric approach of 
\cite{ATV1},\cite{ATV2} was developed it was not even known that this algebra was 
noetherian.  

The classification theory of graded algebras of dimension three, or noncommutative 
projective surfaces, has also progressed substantially in recent years; for a survey of some
of these results see \cite{StVdB}.  In this paper, we study a class of algebras of dimension at least 
three and show that they provide counterexamples to a number of open questions in the 
literature.  In particular, these give new examples of noncommutative surfaces with 
much different behavior than any of those studied previously.  Moreover, the 
construction of these algebras is quite simple and general.

\begin{definition}
\label{const R} $S = \bigoplus_{i \geq 0} S_i$ be a generic Zhang twist (see 
\S\ref{Zhang twists}) of a polynomial ring in $(t+1)$ variables over an algebraically closed field 
$k$ for some $t \geq 2$.  Then let $R$ be the subalgebra of $S$ generated by any generic 
codimension-one subspace of $S_1$.  
\end{definition}
To give a very explicit example of the surface case, for $t = 2$ we may take 
\[
S = k\{x, y, z | xz = pzx, yz = qzy, xy = pq^{-1}yx \} \] for any scalars $p,q \in k$ 
which are algebraically independent over the prime subfield of $k$.  Then let $R$ be 
the subalgebra of $S$ generated over $k$ by any two linearly independent elements $r_1, r_2 
\in S_1$ such that the $k$-span of $r_1$ and  $r_2$ does not contain $x$, $y$, or $z$.

For the rest of this introduction, we assume that all algebras $A = 
\oplus_{i=0}^{\infty} A_i$ are $\mb{N}$-graded and finitely generated by $A_1$ as an 
algebra over $A_0 = k$, where $k$ is an algebraically closed field.  
A \emph{point module} over $A$ is a cyclic $\mb{N}$-graded 
left module $M$, generated in degree 0, such that $\dim_k M_i = 1$ for all $i \geq 
0$.  In case $A$ is commutative, the isomorphism classes of point modules over $A$ are in natural correspondence with the 
closed points of the scheme $\rprojnodash A$.  More generally, for many specific examples of 
noncommutative graded rings one may show explicitly that 
the set of point modules is parameterized by a 
commutative scheme, and the geometry of this scheme often gives important information about 
the ring itself; for example see \cite{ATV1}, \cite{ATV2}.
This is a very useful technique, and so 
it is natural to wonder in what generality the point modules for a ring will form a nice geometric 
object.  Let a $k$-algebra $A$ be called \emph{strongly (left) noetherian} if $A \otimes_k B$ is a 
left noetherian ring for all commutative noetherian $k$-algebras $B$.  A recent theorem of Artin and Zhang
(see Theorem~\ref{pm cor intro} below) shows that the point modules for any strongly noetherian $k$-algebra
are naturally parameterized by a commutative projective scheme over $k$.

The strong noetherian property holds for many standard examples of noncommutative 
rings, including all finitely generated commutative $k$-algebras, all twisted 
homogeneous coordinate rings of projective $k$-schemes, and the AS-regular algebras 
of dimension 3 \cite[Section 4]{ASZ}.  On the other hand, Resco and Small \cite{RSm} have 
given an example of a noetherian finitely generated algebra over a field which is not 
strongly noetherian.  This algebra is not graded, however, nor is the base field 
algebraically closed, and so the example falls outside the paradigm of noncommutative 
projective geometry.  Artin, Small and Zhang have asked if perhaps every 
finitely generated $\mb{N}$-graded noetherian $k$-algebra is strongly noetherian.
We will prove the following theorem which answers this question in the negative.
\begin{theorem} 
(Theorem~\ref{not strongly noetherian}) \label{main result 1} The ring $R$ of 
Definition~\ref{const R} is a connected $\mb{N}$-graded $k$-algebra, finitely 
generated in degree 1, which is noetherian but not strongly noetherian.  \hfill $\Box$
\end{theorem}

We offer two different proofs that $R$ is not strongly noetherian.  
First, we will classify the set of point modules for $R$, from which we can see that 
$R$ fails to satisfy Artin and Zhang's theorem (Theorem~\ref{pm cor intro}).  
For the second proof, we construct an explicit commutative noetherian ring $B$ such 
that $R \otimes_k B$ is not noetherian.  The ring $B$ that works is an infinite 
affine blowup of affine space, which was defined in \cite{ASZ} and is an interesting construction 
in itself.

In \cite{AZ94}, Artin and Zhang describe a categorical approach to noncommutative 
geometry.  Let $A\gr$ be the category of all noetherian graded left $A$-modules, and 
let $A\tors$ be the subcategory of all modules with finite $k$-dimension.  If $A$ is 
commutative, then part of Serre's theorem states that the category $\coh X$ of 
coherent sheaves on the commutative projective scheme $X = \rprojnodash A$ is 
equivalent to the quotient category $A\qgr = A\gr/A\tors$ \cite[Exercise II.5.9]{H}.  
Using this as motivation, for any graded algebra $A$ the \emph{noncommutative 
projective scheme} $A\proj$ is defined to be the pair $(A\qgr, \pi(A))$, where  
$\pi(A)$ is the image of $A$ in $A\qgr$ and plays the role of the structure sheaf.  
One defines cohomology groups for $\mc{M} \in A\qgr$ by setting $\coH^i(\mc{M}) = 
\ext^i_{A\qgr} (\pi(A), \mc{M})$. 

Some homological conditions on the ring $A$ called the $\chi$ conditions arise 
naturally in this approach.  Let $_A k = A/\oplus_{n=1}^{\infty} A_n$; then we say 
that $A$ satisfies $\chi_i$ if $\dim_k \uext^j_A(k,M) < \infty$ for all $0 
\leq j \leq i$ and all $M \in A\gr$.  We say that $A$ satisfies $\chi$ if $A$ 
satisfies $\chi_i$ for all $i \geq 0$.  The most important of these conditions is 
$\chi_1$:  if $A$ satisfies $\chi_1$ then a noncommutative version of Serre's theorem 
holds (see Theorem~\ref{noncomm Serre} below), which shows that the ring $A$ is 
determined in large degree by its associated scheme $A\proj$ together with the 
natural shift functor. 

The $\chi$ conditions hold trivially for commutative rings, but Stafford and Zhang 
\cite[Theorem~2.3]{StZh} gave an example of a noetherian ring $T$ of dimension $2$ 
for which $\chi$ fails.  The ring $T$ fails $\chi_i$ for all $i \geq 0$, and so it 
does not satisfy the noncommutative Serre's theorem.  The algebra $R$ of 
Definition~\ref{const R} is a more interesting example of the failure of $\chi$, 
since $\chi_1$ holds and the noncommutative Serre's theorem is valid for $R$.  
The following 
theorem thus answers a question of Stafford and Zhang from \cite[Section 4]{StZh}.
\begin{theorem} (Theorem~\ref{chi for R})
\label{main result 2} $R$ is a noetherian connected finitely $\mb{N}$-graded 
$k$-algebra, finitely generated in degree 1, for which $\chi_1$ holds but $\chi_i$ 
fails for all $i \geq 2$. \hfill $\Box$
\end{theorem}
One consequence we will draw is that the category $R\qgr$ is necessarily quite different 
from the standard examples of noncommutative schemes, which come from rings satisfying $\chi$.
\begin{theorem}
(Theorem~\ref{R proj not comm}) 
\label{main result 4} The category $R\qgr$ is not equivalent to the category $\coh X$ 
of coherent sheaves on $X$ for any commutative projective scheme $X$.  More 
generally, $R\qgr \not \sim A\qgr$ for any graded ring $A$ which satisfies 
$\chi_2$.  \hfill $\Box$
\end{theorem}

A \emph{maximal order} is the noncommutative analogue of a commutative integrally 
closed domain; see \S \ref{max order section} for the formal definition.  In 
\cite[page 194]{StVdB}, the authors ask whether the $\chi$ conditions perhaps must 
hold for all maximal orders, one reason being that all of the noetherian examples in 
\cite{StZh} for which $\chi$ fails are equivalent orders to maximal orders which do 
satisfy $\chi$.  We show to the contrary the following result.  
\begin{theorem}(Theorem~\ref{max order}, Theorem~\ref{chi for R})
\label{main result 3} $R$ is a connected, finitely $\mb{N}$-graded maximal order for 
which $\chi_i$ fails for $i \geq 2$. \hfill $\Box$
\end{theorem}

For a graded ring $A$, the \emph{global dimension} of $A\qgr$ is the supremum of all 
$n$ such that $\ext^n_{A\qgr}(\mc{M}, \mc{N}) \neq 0$ for some $\mc{M}, \mc{N} \in 
A\qgr$.  Similarly, we define the \emph{cohomological dimension} of $A\proj$ to be 
the supremum of all $n$ such that $\coH^n(\mc{F}) \neq 0$ for some  $\mc{F} \in 
A\qgr$.  Despite the failure of $\chi$ for $R$,  the following theorem shows that 
$R\qgr$ is not too badly behaved.  
\begin{theorem}(Theorem~\ref{cohom dim R})
\label{main result 5} $R\qgr$ has global dimension at most $t +1$ and cohomological 
dimension at most $t$, where $t = \GK(R) -1$.  \hfill $\Box$
\end{theorem}

The starting point for our study was the work of David Jordan on rings generated by 
two Eulerian derivatives \cite{Jordan}.  These rings are more or less the algebras 
$R$ of Definition~\ref{const R} for $t = 2$.  Jordan had classified the point modules 
for such algebras and showed that the strong noetherian property must fail, but did 
not show these algebras were noetherian.  Our results imply in 
fact that rings generated by a generic finite set of Eulerian derivatives are 
noetherian, answering a question in \cite{Jordan}---see \S \ref{eul der} for more 
details.  

The results in this paper form part of the author's PhD thesis \cite{Rothesis}, and 
in some cases extra details may be found there.  In collaboration with J.T. Stafford 
and D.S. Keeler we have recently developed a geometric approach to the study of the rings
$R$ which shows that they may be thought of as a kind of peculiar noncommutative blowing up 
of projective space; details will be given in a further paper.  The author would like to thank
Jessica Sidman, Daniel Chan, and Melvin Hochster for helpful discussions, and David Jordan, whose 
work in \cite{Jordan} motivated the study of these algebras.  Toby Stafford deserves special 
thanks for his invaluable advice and support.

%_____________________________________________________________
\section{Basic definitions} \label{definitions} In this section, 
we fix some terminology concerning noncommutative graded rings.  
The reader may wish to skim this section and refer back to it later 
when necessary. 

We make the convention throughout that $0$ is a natural number, so that $\mb{N} = 
\mb{Z}_{\geq 0}$. The rings $A$ we study in this paper are all $\mb{N}$-graded 
algebras $A = \oplus_{i = 0}^{\infty} A_i$ over an algebraically closed field $k$.  
We assume throughout that all algebras $A$ are \emph{connected}, that is that $A_0 = 
k$, and finitely generated as an algebra by $A_1$.  Let $A\Gr$ be the abelian 
category whose objects are the $\mb{Z}$-graded left $A$-modules $M = \oplus_{i = 
-\infty}^{\infty} M_i$, and where the morphisms $\caphom_{A\Gr}(M,N)$ are the module 
homomorphisms $\phi$ satisfying $\phi(M_n) \subseteq N_n$ for all $n$.  Let $A\gr$ be 
the full subcategory of $A\Gr$ consisting of the noetherian objects.  For $M \in 
A\Gr$ and $n \in \mb{Z}$, the shift of $M$ by $n$, denoted $M[n]$, is the module with 
the same ungraded module structure as $M$ but with the grading shifted so that 
$(M[n])_m = M_{n+m}$.  Then for $M, N \in A\Gr$ we may define
\[
\uhom_A(M,N) = \oplus_{i = - \infty}^{\infty} \caphom_{A\Gr}(M, N[i]).
\]
The group $\uhom(M,N)$ is a $\mb{Z}$-graded 
vector space and we also write $\uhom(M,N)_i$ for 
the $i$th graded piece $\caphom(M,N[i])$.  Under mild hypotheses, for example if $M$ 
is finitely generated, the group $\uhom(M,N)$ may be identified with the set of 
ungraded module homomorphisms from $M$ to $N$.  It is a standard result that the 
category $A\Gr$ has enough injectives and so we may define right derived functors 
$\ext^i_{A\Gr}(M, -)$ of $\caphom_{A\Gr}(M, -)$  for any $M$.  The definition of $\uhom$ 
generalizes to
\[
\uext^i_A(M,N) =\oplus_{i = - \infty}^{\infty} \ext^i_{A\Gr}(M, N[i]).
\]
See \cite[Section 3]{AZ94} for a discussion of the basic properties of $\uext$.   

For a module $M \in A\Gr$, a \emph{tail} of $M$ is any submodule of the form $M_{\geq 
n} = \oplus_{i = n}^{\infty} M_i$.  A \emph{subfactor} of $M$ is any module of the 
form $N/N'$ for graded submodules $N' \subseteq N$ of $M$.  For the purposes of this 
paper, $M \in A\Gr$ is called \emph{torsion} if for all $m \in M$ there exists some 
$n \geq 0$ such that $(A_{\geq n}) m = 0$.  Note that if $M \in A\gr$, then $M$ is 
torsion if and only if $\dim_k M < \infty$.  We say that $M \in A\Gr$ is \emph{left 
bounded} if $M_i = 0$ for $i \ll 0$, and \emph{right bounded} if $M_i = 0$ for $i \gg 
0$. $M$ is \emph{bounded} if it is both left and right bounded.  
A \emph{(finite) filtration} of $M \in A\Gr$ is a sequence of graded submodules $0 = 
M_0 \subseteq M_1 \subseteq \dots \subseteq M_n = M$; we call the modules 
$M_i/M_{i-1}$ for $0 < i \leq n$ the \emph{factors} of the filtration. 

A \emph{point module} over $A$ is a graded module $M$ such that $M$ is cyclic, 
generated in degree $0$, and $\dim_k M_n = 1$ for all $n \geq 0$.  Note that a tail 
of a point module is a shift of some other point module.  A $\emph{point ideal}$
is a left ideal $I$ of $A$ such that $A/I$ is a point module, or equivalently
such that $\dim_k I_n = \dim_k A_n -1$ for all $n \geq 0$.  Since $A$ is generated 
in degree $1$, the point ideals of $A$ are in one-to-one correspondence with
isomorphism classes of point modules over $A$. 

We will use Gelfand-Kirillov dimension, or GK-dimension for short, as our dimension 
function for modules; the basic reference for its properties is \cite{KL}.  Given $M 
\in A\gr$, the \emph{Hilbert function} of $M$ is the function $H(n) = \dim_k M_n$ for 
$n \in \mb{Z}$.  If $A$ is a finitely generated $k$-algebra and $_A M$ is a finitely generated
module, then 
$\GK(M)$ depends only on the Hilbert function of $M$ \cite[6.1]{KL}.  In particular, 
if $M$ has a \emph{Hilbert polynomial}, that is $\dim_k M_n = f(n)$ for $n \gg 0$ and 
some polynomial $f \in \mb{Q}[n]$, then $\GK(M) = \deg f +1$ (with the convention 
$\deg (0) = -1$).  We say $M \in A\Gr$ is \emph{(graded) critical} if $\GK(M/N) < 
\GK(M)$ for all nonzero graded submodules $N$ of $M$. 
If $A$ is an 
$\mb{N}$-graded noetherian ring, then the GK-dimension for $A$-modules is 
\emph{exact}: in other words, given any exact sequence $0 \ra M' \ra M \ra M'' \ra 0$ 
in $A\gr$, one has $\GK(M) = \max (\GK(M'), \GK(M''))$ \cite[4.9]{MR}.
 
%_______________________________________________________________________________
\section{Zhang twists}
\label{Zhang twists} 
Let $k$ be an algebraically closed field.  
Fix from now on a commutative polynomial ring $U = k[x_0, \dots, x_t]$ in $t+1$ indeterminates, 
graded as usual with $\deg x_i = 1$ for all $i$, and the corresponding projective 
space $\operatorname{Proj} U = \mb{P}^t$.  By a \emph{point} of $\mb{P}^t$ we always 
mean a closed point.  The main results of this paper require that $t \geq 2$, so we assume this throughout.

Let the symbol $\regmult$ indicate the multiplication operation in $U$.  For any 
graded automorphism $\phi$ of $U$ we may define a new graded ring $(S, \twistmult)$, 
where $S$ has the same underlying vector space as $U$ and $f \twistmult g = \phi^n(f) 
\regmult g$ for $f \in U_m$, $g \in U_n$.  This is just a special case of the (left) 
twisting construction studied by Zhang in \cite{Zhang}.  

Let $\m_d$ stand for the ideal in $U$ of a point $d \in \mb{P}^t$.  For a homogeneous element $f \in U$, we
use the notation $f \in \m_d$ and $f(d) = 0$ interchangeably to indicate that $f$ vanishes at $d$.  
Corresponding to the
automorphism $\phi$ of $U$ is an automorphism $\varphi$ of 
$\mb{P}^t$ which satisfies $\mf{m}_{\varphi(d)} = \phi^{-1}(\mf{m}_d)$ for all $d \in 
\mb{P}^t$.  Automorphisms $\phi_1, \phi_2$ of $U$ give the same automorphism 
$\varphi$ of $\mb{P}^t$ if and only if $\phi_1 = a \phi_2$ for some $a \in 
k^{\times}$ \cite[II.7.1.1]{H}.  Moreover, automorphisms of $U$ which are scalar 
multiples give isomorphic twisted algebras $S$ \cite[5.13]{Zhang}.  Thus a particular 
twist $S$ of $U$ is determined up to isomorphism by $\varphi$ and 
we write $S = S(\varphi)$.  We remark that 
an alternative description of $S$ may be given using 
twisted homogeneous coordinate rings \cite[Section 3]{StVdB}.  In this language, 
$S = B(\mb{P}^t, \mc{O}(1), \varphi^{-1})$ where $\mc{O}(1)$ is the twisting sheaf
of Serre on $\mb{P}^t$.    

Since $S$ and its subalgebras are our main interest in this paper, our notational convention from now 
on (except in the appendix) will be to let juxtaposition indicate multiplication in 
$S$ and to use $\regmult$ when multiplication in $U$ is intended.  
However, we let exponents retain their old commutative meaning, so 
that if $\m_d$ is a point ideal of $U$ then $\m_d^2$ is short for $\m_d \circ \m_d$, the 
ideal of polynomials vanishing twice at $d$.  

Many important properties pass from a graded ring to any Zhang twist.  In particular, 
it is easy to see that $S$ is a noetherian domain of GK dimension $t+1$, since these 
properties are obvious for $U$ \cite[Propositions~5.1,5.2,5.7]{Zhang}.  Suppose that 
$M \in U\Gr$, and let $\regmult$ indicate the action of $U$ on $M$. Then $M$ obtains 
an $S$-structure using the rule $sx = \phi^n(s) \regmult x$ for $s \in S_m$, $x \in 
M_n$.  This defines a functor $\theta: U\Gr \ra S\Gr$ which is an 
equivalence of categories \cite[3.1]{Zhang}.  For any graded ideal $I$ of $U$, 
$\theta(I)$ is a graded left ideal of $S$.  We often simply identify the underlying $k$-spaces 
of these ideals and call both $I$.  

Since the equivalence of categories preserves Hilbert functions and the property of being cyclic,
it is clear that 
the point modules over $S$ are the modules of the form $\theta(U/\m_d) = 
S/\m_d$ for $d \in \mb{P}^t$.  We will use the following notation:
\begin{notation}
\label{point mod} Given a point $d \in \mb{P}^t$, let $\pointmod{d}$ be the left point 
module $\theta(U/\m_d) = S/\m_d$ of $S$.
\end{notation}

We may also describe point modules over $S$ by their \emph{point sequences}.  If $M$ 
is a point module over $S$, then the annihilator of $M_n$ in $S_1 = U_1$ is some 
codimension 1 subspace which corresponds to a point $d_n \in \mb{P}^t$.  The point sequence of $M$ 
is defined to be the sequence $(d_0, d_1, d_2, \dots) $ of points of $\mb{P}^t$.  
Clearly two $S$-point modules are isomorphic if and only if they have the same point 
sequence.
\begin{lemma}
\label{truncating point modules}  Let $d$ be an arbitrary point of $\mb{P}^t$. 
\begin{enumerate}
\item \label{truncating point modules 1} 
$\pointmod{d}$ has point sequence $(d, \varphi(d), \varphi^2(d), \dots)$.  
\item \label{truncating point modules 2}
$(\pointmod{d})_{\geq n} \cong \pointmod{\varphi^n(d)}[-n]$ as $S$-modules.
\end{enumerate} 
\end{lemma}
\begin{proof}
(1) By definition $P(d) = S/\m_d$.  If $f \in S_1$, then $f S_n = \phi^n(f) \regmult 
U_n \subseteq \mf{m}_d$ if and only if $\phi^n(f) \in \mf{m}_d$, in other words $f 
\in \phi^{-n}(\mf{m}_d) = \mf{m}_{\varphi^n(d)}$.

(2) By part (1), $(\pointmod{d})_{\geq n}$ is the shift by $(-n)$ of the point module 
with point sequence $(\varphi^n(d), \varphi^{n+1}(d), \dots)$.   
\end{proof}

Finally, we record the following simple facts which we shall use frequently.
\begin{lemma}
\label{cyclic critical filter} Let $M \in S\gr$.
\begin{enumerate} 
\item \label{cyclic critical filter 1}
If $M$ is cyclic $1$-critical, then $M \cong \pointmod{d}[i]$ for some $d \in 
\mb{P}^t$ and $i \in \mb{Z}$.  
\item \label{cyclic critical filter 2} 
$M$ has a finite filtration with factors which are graded cyclic critical $S$-modules.
\end{enumerate}
\end{lemma}
\begin{proof}
(1) The equivalence of categories $U\Gr \sim S\Gr$ preserves the GK-dimension of 
finitely generated modules, since it preserves Hilbert functions, and so it also 
preserves the property of being $\GK$-critical.  It is standard that the cyclic 
$1$-critical graded $U$-modules are just the point modules over $U$ and their 
shifts.  Under the equivalence of categories, the corresponding $S$-modules are the 
$S$-point modules and their shifts.   

(2) Each module $N \in U\gr$ has a finite filtration composed of graded cyclic 
critical $U$-modules, so the same holds for $S$-modules by the equivalence of 
categories.  
\end{proof}

%_______________________________________________________________________________
\section{The algebras $R(\varphi,c)$}

Let $S = S(\varphi)$ for some $\varphi \in \aut \mb{P}^t$.  For any codimension-1 
subspace $V$ of $S_1 = U_1$, we let $R = k \langle V \rangle \subseteq S$ be the 
subalgebra of $S$ generated by $V$.  The vector subspace $V$ of $U_1$ corresponds to 
a unique point $c \in \mb{P}^t$.   Then $R$ is determined up to isomorphism by the geometric 
data $(\varphi,c)$ and we write $R = R(\varphi,c)$.  We emphasize again that
we always assume that $t \geq 2$ from now on; for smaller $t$ the ring $R$ 
is not very interesting.

We shall see that the basic properties of $R(\varphi,c)$ depend closely on properties
of the iterates of the point $c$ under $\varphi$.  It is convenient to let $c_i = 
\varphi^{-i}(c)$ for all $i \in \mb{Z}$.  Then the ideal of the point $c_i$ is 
$\phi^i(\mf{m}_c)$.  In case $c$ has finite order under $\varphi$, that is 
$\varphi^n(c) = c$ for some $n > 0$, the algebra $R(\varphi,c)$ behaves quite 
differently from the case where $c$ has infinite order.  For example, if $\varphi^n 
=1$ for some $n \geq 0$ then it is easy to see that $S$ and hence $R$ are $PI$ 
rings.  The finite order case turns out to have none of the interesting properties of 
the infinite order case, 
and so in this paper we will only 
consider the case where $c$ has infinite order under $\varphi$.
\begin{hypothesis}
\label{distinct} Assume that $(\varphi,c) \in (\aut \mb{P}^t) \times 
\mb{P}^t$ is given such that $c$ has infinite order under $\varphi$, or equivalently 
the points $\{c_i\}_{i \in \mb{Z}}$ are distinct.
\end{hypothesis}

We note some relationships among the various $R(\varphi,c)$.  In particular, part (1)
of the next lemma will allow
us to transfer our left sided results to the right. 
\begin{lemma}
\label{switch sides} Let $(\varphi,c) \in (\aut \mb{P}^t) \times \mb{P}^t$, and let 
$\psi$ be any automorphism of $\mb{P}^t$.  Then 
\begin{enumerate}
\item \label{switch sides 1}
$R(\varphi, c)^{op} \cong R(\varphi^{-1}, \varphi(c))$. 
\item \label{switch sides 2} 
$R(\varphi, c) \cong R(\psi\varphi\psi^{-1}, \psi(c))$.
\end{enumerate}
\end{lemma}
\begin{proof} 
(1) Set $S = S(\varphi)$ and $S' = S(\varphi^{-1})$, identifying the underlying 
spaces of each with that of $U$.  Let $\phi$ be an automorphism of $U$ corresponding 
to $\varphi$.  Then it is straightforward to check that the vector space map defined 
on the graded pieces of $U$ by sending $f \in U_m$ to $\phi^{-m}(f) \in U_m$ is a 
graded algebra isomorphism from $S^{op}$ to $S'$.  The isomorphism maps $(\m_c)_1$ to 
$(\m_{\varphi(c)})_1$ and so it restricts to an isomorphism $R(\varphi, c)^{op} \cong 
R(\varphi^{-1}, \varphi(c))$. 

(2) Similarly, let $\sigma$ be an automorphism of $U$ corresponding to $\psi$.  It is 
easy to check that the vector space map of $U$ defined by $f \mapsto \sigma^{-1}(f)$ 
is an isomorphism of $S(\varphi)$ onto $S(\psi\varphi\psi^{-1})$ which maps 
$(\m_c)_1$ to $\sigma^{-1}(\m_c)_1 = (\m_{\psi(c)})_1$, and so restricts to an 
isomorphism $R(\varphi,c) \cong R(\psi \varphi \psi^{-1}, \psi(c))$.
\end{proof}

In the next theorem we will prove an important characterization of the elements of $R 
= R(\varphi, c)$ which is foundational for all that follows.  First we need the 
following lemma; the proofs of it and several other technical commutative results which appear later in 
the paper may be found in the appendix.  

\begin{lemma}(Lemma~\ref{products versus intersections app})
\label{products versus intersections} Let $\m_1, \m_2, \dots, \m_n$ be the ideals of $U$ 
corresponding to distinct points $d_1, \dots d_n$ in $\mb{P}^t$.  Then $(\m_1 \circ 
\m_2 \circ \dots \circ \m_n)_{\geq n} = (\bigcap_{i=1}^n \m_i)_{\geq n}$. \hfill 
$\Box$
\end{lemma}

\begin{theorem}
\label{char of R} Let $R = R(\varphi, c)$.  Then for all $n \geq 0$,
\[
R_n = \{f \in U_n \mid f(c_i) = 0\ \text{for}\ 0 \leq i \leq n-1 \}.
\]
\end{theorem}

\begin{proof}  By definition $R = k \langle V \rangle \subseteq S$, where $V =
(\mf{m}_c)_1$ considered as a subset of $U$.  For $n = 0$ the statement of the
theorem is $R_0 = U_0 = k$, which is clearly correct, so assume that $n \geq 1$.  Then
\[ R_n = V^n = \phi^{n-1}(V) \regmult
\phi^{n-2}(V) \regmult \dots \regmult \phi(V) \regmult V. \] Now $\phi^i(V) =
(\mf{m}_{c_i})_1$, and the points $c_i$ are distinct by Hypothesis~\ref{distinct}.  
Thus by Lemma~\ref{products versus intersections} we get that 
\begin{gather*}
R_n = (\mf{m}_{c_{n-1}})_1 \regmult \dots \regmult (\mf{m}_{c_1})_1 \regmult
(\mf{m}_{c_0})_1 = [(\mf{m}_{c_{n-1}}) \regmult \dots \regmult (\mf{m}_{c_1})
\regmult (\mf{m}_{c_0})]_n \\
= [(\mf{m}_{c_{n-1}}) \cap \dots \cap (\mf{m}_{c_1}) \cap (\mf{m}_{c_0})]_n.
\end{gather*}
The statement of the theorem in degree $n$ follows.
\end{proof}

The theorem has a number of easy consequences.  
The first will be a simple calculation of the Hilbert function of $R$, which depends on the following fundamental
commutative result which is proved in the appendix.
\begin{lemma} (Lemma~\ref{Hilbert of point ideal app})
\label{Hilbert of point ideal} let $d_1, d_2,\dots d_n$ be distinct points in
$\mb{P}^t$, for some $n \geq 1$, and let $\m_1, \m_2, \dots, \m_n$ be 
the corresponding graded ideals of $U$.  Let $e_i > 0$ for all $1 
\leq i \leq n$.  Set $J = \bigcap_{i=1}^n \m_i^{e_i}$. Then $\dim_k J_m = 
\binom{m+t}{t} - 
\sum_i \binom{e_i+t-1}{t}$ for all $m \geq (\sum e_i) -1$. \\
In particular, if $J = 
\bigcap_{i=1}^n \m_i$ then $\dim_k J_m = \binom{m+t}{t} - n$ for $m \geq n-1$. 
\hfill $\Box$
\end{lemma}

\begin{lemma}
\label{Hilbert R}
Let $R = R(\varphi,c)$.  Then $\dim_k R_n = \binom{n+t}{t} - n$ for all $n \geq 0$.  
In particular, $\GK(R) = t+1$.
\end{lemma}
\begin{proof} 
The Hilbert function of $R$ follows from Theorem~\ref{char of R} and 
Lemma~\ref{Hilbert of point ideal}.  Since we always assume that $t \geq 2$, it is 
clear that the Hilbert polynomial of $R$ has degree $t$ and so $\GK(R) = t+1$.
\end{proof}

\begin{lemma}
\label{quotient rings} 
\label{essential} 
The rings $R = R(\varphi,c)$ and $S = S(\varphi)$ have the same graded quotient 
ring $D$ and Goldie quotient ring $Q$.  The inclusion $R \hookrightarrow 
S$ is a essential extension of left (or right) $R$-modules.
\end{lemma}
\begin{proof} Since both $R$ and $S$ are domains of finite GK-dimension, they both have
graded quotient rings and Goldie quotient rings \cite[C.I.1.6]{NV}, \cite[4.12]{KL}.
Clearly the graded quotient ring $D'$ of $R$ is contained in the graded quotient ring 
$D$ of $S$.  Since we assume always that $t \geq 2$, we may choose a nonzero 
polynomial $g \in S_1$ with \mbox{$g \in \mf{m}_{c_0} \cap \mf{m}_{c_1}$.}  
Then Theorem~\ref{char of R} implies that $g \in R_1$ and $S_1 g \subseteq R_2$.  Thus 
$S_1 \subseteq R_2 (R_1)^{-1} \subseteq D'$ and consequently $D' = D$. Then $Q$, the 
Goldie quotient ring of the domain $D$, is also the Goldie quotient ring for both $R$ 
and $S$.  The last statement of the proposition is now clear.
\end{proof}

%________________________________________________________________________________
\section{The noetherian property for $R$}
\label{R noeth}

Let $S = S(\varphi)$ and $R = R(\varphi, c)$.  In this section we will 
characterize those choices of $\varphi$ and $c$ satisfying 
Hypothesis~\ref{distinct} for which the ring $R$ is noetherian.  To do this,
we will first analyze the structure of the factor module $_R (S/R)$ in detail, and then use
this information to understand contractions and extensions of left ideals between $R$ and $S$.  

The following notation will be convenient in this section.
\begin{notation}
\label{set not} 
\begin{enumerate}
\item \label{set not 1}
$A_n = \{0,1, \dots n-1\}$ for $n > 0$ and $A_n = \emptyset$ for $n \leq 0$. 
\item \label{set not 2}
For $B \subseteq \mb{Z}$, set $B + m = \{b + m \mid b \in B \}$.
\end{enumerate}
\end{notation}

\begin{definition}
\label{special modules} Let $B \subseteq \mb{N}$.  
We define a left $R$-module $\Tmod{B} \subseteq S$ by specifying its
graded pieces as follows:
\[
(\Tmod{B})_n = \{f \in S_n\ \mid f(c_i) = 0\ \text{for}\ i \in A_n \setminus B \}
\]
We then define the left $R$-module $\specmod{B} = \Tmod{B}/R \subseteq (S/R)$.
\end{definition}

We should check that $\Tmod{B}$ really is closed under left multiplication by $R$.  
If $g \in R_m$ and $f \in (\Tmod{B})_n$, then $gf = \phi^{n}(g) \regmult  f$.  Now 
$g(c_i) = 0$ for $i \in A_m$ by Theorem~\ref{char of R} and $f(c_i) = 0$ for $i \in 
A_n \setminus B$ by definition.  Thus $[\phi^{n}(g) \regmult  f](c_i) = 0$ for $i \in 
(A_m + n) \cup (A_n \setminus B) \supseteq (A_{n+m} \setminus B$), and so $gf \in 
(\Tmod{B})_{n+m}$ as required.  Also, by Theorem~\ref{char of R} the extreme cases 
are $R = \Tmod{\emptyset}$ and $S = \Tmod{\mb{N}}$.  In particular, $R \subseteq 
\Tmod{B}$ always holds, and so $\specmod{B}$ is well defined.  

\begin{lemma}
\label{Hilbert M_B} The Hilbert function of $\specmod{B}$ is given by
\[
\dim_k (\specmod{B})_n = |A_n \cap B|.
\]
\end{lemma}
\begin{proof} Immediate from Lemma~\ref{Hilbert of point ideal}.
\end{proof}

In the special case of Definition~\ref{special modules} where $B$ is a singleton set, 
$\specmod{B}$ is just a shifted $R$-point module, as follows.
\begin{lemma}
\label{one point} Let $j \in \mb{N}$.  Then $M = \specmod{\{j\}}$ is an $R$-point
module shifted by $j+1$.  In fact, $M \cong {} _R \pointmod{c_{-1}}[-j-1]$.
\end{lemma}

\begin{proof}
By Lemma~\ref{Hilbert M_B} the Hilbert function of $M$ is
\[
\dim_k M_n =
\begin{cases}
0 & 0 \leq n \leq j \\
1 & j+1 \leq n
\end{cases}
\]
so that $M$ does have the Hilbert function of a point module shifted by $j+1$.  For 
convenience of notation set $m = j+1$, and let us calculate $\ann_R (M_m)$.  Now $f 
M_m = 0$ for $f \in R_n$ if and only if $f (\Tmod{\{j\}})_m \subseteq R_{m+n}$.  
Since $M_m \neq 0$, we may choose $g \in (\Tmod{\{j\}})_m$ such that $g \not \in R$; 
then $g(c_j) \not = 0$.  Also, because $\dim_k M_m = 1$ it is clear that 
$(\Tmod{\{j\}})_m = R_m + kg$, and so $f (\Tmod{\{j\}})_m \subseteq R$ if and only if 
$fg \in R$.  Now $fg = \phi^m(f) \regmult g$, and so by Theorem~\ref{char of R} we 
have that  $fg \in R$ if and only if $\phi^m(f)(c_j) = 0$, equivalently $f(c_{-1}) = 
0$, since $m = j+1$.  In conclusion, $\ann_R(M_m) = \mf{m}_{c_{-1}} \cap R$.   

Thus we have an injection of $R$-modules given by right multiplication by $g$:
\[
\psi: (R/(\m_{c_{-1}} \cap R))[-m] \overset{g}{\lra} \Tmod{\{j\}}/R = M.
\]
By Lemma~\ref{Hilbert of point ideal}, $R/(\m_{c_{-1}} \cap R)$ has the Hilbert 
function of a point module and so both sides have the same Hilbert function.  Thus 
$\psi$ is actually an isomorphism.  In particular, $M$ is cyclic and so is a shifted 
$R$-point module. 

We also have the injection $R/(\mf{m}_{c_{-1}} \cap R) \ra S/(\mf{m}_{c_{-1}}) = 
\pointmod{c_{-1}}$, and since both sides have the Hilbert function of a point 
module this is also an isomorphism of $R$-modules.  So $M \cong {}_R 
(\pointmod{c_{-1}})[-j-1]$.
\end{proof}

We may now understand the structure of $_R (S/R)$ completely.
\begin{proposition}
\label{sum of points}
The modules $\{\specmod{\{j\}}\}_{j \in \mb{N}}$
are independent submodules of $S/R$. Also, for $B \subseteq \mb{N}$,
\[ \specmod{B} =
\bigoplus_{j \in B} \specmod{\{j\}}.
\]
\end{proposition}

\begin{proof} 
We first show the independence of the $\specmod{\{j\}}$.  It is enough to work with 
homogeneous elements; fix $n \geq 0$ and let $\sum_{j \in \mb{N}} f_j = 0$ for some 
$f_j \in (\specmod{\{j\}})_n$.  Let $f_j = g_j + R$ for some elements $g_j \in 
(\Tmod{\{j\}})_n \subseteq S_n$.  Thus $\sum g_j \in R_n$.  Suppose that some $f_j 
\not = 0$, and let $k = \min\{j | f_j \not = 0\}.$   Now $(\specmod{\{j\}})_{\leq j} 
= 0$ for all $j$ by Lemma~\ref{Hilbert M_B}, and so we must have $k \leq n -1$. Then 
$k \in A_n \setminus \{j\}$ for any $j \not = k$, so that $g_j(c_k) = 0$ for all $j 
\not = k$, by the definition of $\Tmod{\{j\}}$.  Since $\sum g_j \in R_n$, 
Theorem~\ref{char of R} implies that $(\sum g_j)(c_k) = 0$ also holds and so 
$g_k(c_k) = 0$. But then $g_k \in (\Tmod{\{k\}})_n \cap ({\mf{m}_{c_k}}) = R_n$, and 
thus $f_k = 0$, a contradiction.  We conclude that all $f_j =0$, and so the 
$\specmod{\{j\}}$ are independent.

For the second statement of the proposition, by Lemma~\ref{Hilbert M_B} the Hilbert 
function of $\specmod{B}$ is $\dim_k (\specmod{B})_n = |A_n \cap B|$, while the 
Hilbert function of $\bigoplus_{j \in B} \specmod{\{j\}}$ is
\[
\dim_k \big[ \bigoplus_{j \in B} \specmod{\{j\}} \big]_n = \#\{j \in B | j \leq n-1
\} =|A_n \cap B|.
\]
Thus the Hilbert functions are the same on both sides of our claimed equality. Since
$\sum_{j \in B} \specmod{\{j\}} \subseteq \specmod{B}$ is clear and we know that the
$\specmod{\{j\}}$ are independent by the first part of the proposition, the equality
follows.
\end{proof}

\begin{corollary}
\label{when noetherian} \label{structure of S/R} 
\begin{enumerate} 
\item \label{structure of S/R 1} 
Given $B \subseteq \mb{N}$, $\specmod{B}$ is a noetherian $R$-module if and only if 
the set $B$ has finite cardinality. 
\item \label{structure of S/R 2} 
$\D _R (S/R) \cong \bigoplus_{j=0}^{\infty} {} _R \pointmod{c_{-1}}[-j-1]$.  In 
particular, $_R (S/R)$ is not finitely generated.
\end{enumerate}
\end{corollary}
\begin{proof} (1) is clear since a point module is noetherian.  (2) follows by taking $B =
\mb{N}$ in the proposition and using also Lemma~\ref{one point}.
\end{proof}

Next, we analyze the
noetherian property for some special types of $R$-modules which may be realized as subfactors of $S/R$.
\begin{proposition}
\label{1st type noetherian} For $f \in R_n$, let $N = (Sf \cap R)/{Rf} \in R\Gr$.  
Set $D = \{ i \in \mb{N} \mid f(c_i) =
0 \}$ and $B = (D-n) \cap \mb{N}$.
 Then
\begin{enumerate}
\item \label{1st type noetherian 1} 
$N \cong \specmod{B}[-n]$ 
\item \label{1st type noetherian 2} 
$N$ is noetherian if and only if $|D| < \infty$. 
\end{enumerate} 
\end{proposition}
\begin{proof} First, if we set $T = \{g \in S \mid gf \in R \}$ and $M = T/R$, then $N \cong
M[-n]$.  So it is enough for (1) to show that $T = \Tmod{B}$.

Let $g \in S_m$ be arbitrary.  Note that $A_n \subseteq D$ since $f \in R_n$. Then
\begin{alignat*}{2}
 & gf  = \phi^n(g) \circ f \in R && \qquad  \\            
\Llra\ & [\phi^n(g) \circ f](c_i) = 0 &&\ \text{for all}\ i \in A_{n+m}  \\
\Llra\ & \phi^n(g)(c_i) = 0 &&\ \text{for all}\ i \in A_{n+m} \setminus D\\
\Llra\ & g(c_i) = 0 &&\ \text{for all}\ i \in (A_{n+m} \setminus D) - n\\
\Llra\ & g(c_i) = 0 &&\ \text{for all}\ i \in A_m \setminus (D-n)\ \ 
 \text{(since $A_n \subseteq D$)} \\
\Llra\ & g \in \Tmod{B} &&\ \text{by Definition~\ref{special modules}}.
\end{alignat*}
Thus $T = \Tmod{B}$ and (1) holds.

For (2), note that $D$ has finite cardinality if and only if $B$ does, and apply
Corollary~\ref{when noetherian}(1).
\end{proof}

\begin{proposition}
\label{2nd type noetherian} For $f \in R_n$, let $M = S/{(R + Sf)} \in R\Gr$.  Set $D 
= \{ i \in \mb{N} \mid f(c_i) = 0 \}$.  Then 
\begin{enumerate}
\item \label{2nd type noetherian 1}
$M \cong \specmod{D}$ 
\item \label{2nd type noetherian 2} 
$M$ is noetherian if and only if $|D| < \infty$.
\end{enumerate}
\end{proposition}
\begin{proof} 
Set $B = {\mb{N} \setminus D}$.  We will show that $R + Sf = \Tmod{B}$.  Then we will 
have that $S/(R + Sf) = S/\Tmod{B} \cong (M^{\mb{N}}/M^B) \cong \specmod{D}$ by 
Proposition~\ref{sum of points}.
 
Suppose that $h \in (R+Sf)_m$; then $h = g_1 + g_2 f = g_1 + \phi^n(g_2) \circ f$ for 
some $g_1 \in R_m$ and $g_2 \in S_{m-n}$.  Now $g_1(c_i) = 0$ for $i \in A_m$, and 
$f(c_i)= 0$ for $i \in D$, so that $h(c_i) = 0$ for $i \in A_m \cap D$.  So we have 
$(R+Sf) \subseteq \Tmod{B}$.  Note that $(R + Sf)_m = R_m = (\Tmod{B})_m$ for $m < 
n$. 

Now, we can calculate the Hilbert function of $R + Sf$.  By Propositions~\ref{1st 
type noetherian}(1) and Lemma~\ref{Hilbert M_B}, we have that $\dim_k ((Sf \cap 
R)/Rf)_m = |A_{m-n} \cap (D-n)| = |A_m \cap D| - n$ for $m \geq n$, since $A_n 
\subseteq D$.  Since the Hilbert functions of $R$, $Sf$, and $Rf$ are all known, one 
may calculate that $\dim_k (R +Sf)_m = \binom{m+2}{2} - |A_m \cap D|$ for $m \geq n$, 
which is equal to $\dim_k (\Tmod{B})_m$.  Thus $R + Sf = \Tmod{B}$.
 
Part (2) is then immediate from Corollary~\ref{when noetherian}(1).
\end{proof}

Given a left ideal $I$ of $R$, we may extend to a left ideal $SI$ of $S$, and then 
contract back down to get the left ideal $SI \cap R$ of $R$.  The factor $(SI \cap 
R)/I$ is built up out of the 2 types of modules we considered in 
Propositions~\ref{1st type noetherian} and \ref{2nd type noetherian}. 
\begin{lemma}
\label{filtration by types}
 Let $I$ be a finitely generated nonzero graded left ideal of $R$,
and set $M = (SI \cap R)/I$.  Then $M$ has a finite filtration $0 = M_0 \subseteq M_1
\subseteq \dots \subseteq M_m = M$ such that each factor $M_{i+1}/M_i$ is isomorphic with
shift to a subfactor of either $(Ss_i \cap R)/Rs_i$ or $S/{(R + Ss_i)}$ for 
some nonzero homogeneous $s_i \in R$.
\end{lemma}
\begin{proof} Let $I = \sum_{i=1}^{n} R r_i$ for some homogeneous $r_i \in R$.  
If $n =1$ the result is obvious, so assume that $n \geq 2$.

Set $J = \sum_{i=1}^{n-1} R r_i$.  By induction on $n$, $(SJ \cap R)/J$ and hence also
its surjective image $(SJ \cap R) + I/I$ have filtrations of the required type.  It 
is enough then to show that 
\[
N = (SI \cap R)/((SJ \cap R) + I) = (SI \cap R)/((SJ + Rr_n) \cap R)
\]
has the required filtration. But $N$ injects into $L = SI/(SJ + Rr_n)$.  Now $R$ is 
an Ore domain by Lemma~\ref{quotient rings}, so we may choose a homogeneous element 
$0 \neq r \in R$ such 
that $rr_n \in J$.  Then $L$ is a surjective image (with shift) of $S/(R + Sr)$, so 
$N$ is a shift of a subfactor of $S/(R + Sr)$.
\end{proof}

In certain circumstances the noetherian property passes to subrings.  The following 
lemma is just a slight variant of a number of similar results in the 
literature (for example, see \cite[Lemma~4.2]{ASZ}).

\begin{lemma}
\label{noetherian passes down} Let $A \hookrightarrow B$ be any extension of
$\mb{N}$-graded rings. Suppose that $B$ is left noetherian, and that $(BI \cap A)/I$ 
is a noetherian left $A$-module for all finitely generated homogeneous left ideals 
$I$ of $A$. Then $A$ is left noetherian.
\end{lemma}

\begin{proof} It is enough to prove that $A$ is graded left noetherian, that is that all
homogeneous left ideals are finitely generated.  Let $I$ be a homogeneous left ideal of
$A$.  Then $BI$ is a homogeneous left ideal of $B$, which is finitely generated since
$B$ is noetherian, and so we may pick a finite set of homogeneous generators $r_1, r_2,
\dots r_n \in I$ such that  $BI = \sum_{i=1}^{n} B r_i$.  Let $J = \sum_{i=1}^{n} A
r_i$.  Then $BI = BJ$, and since $J$ is finitely generated over $A$ we may apply the
hypothesis to conclude that ${(BJ \cap A)}/{J} = {(BI \cap A)}/{J}$ is a noetherian
$A$-module.  The submodule $I/J$ of ${(BI \cap A)}/{J}$ is then noetherian over $A$, in
particular finitely generated over $A$.  Finally, since $J$ is finitely generated over
$A$, so is $I$.
\end{proof}

We note the definition of an unusual geometric condition on a set of points of a 
variety, which appeared in \cite[p. 582]{ASZ}.
\begin{definition}
\label{crit dense} \label{criticallydense-index} Let $\mc{C}$ be an infinite set of 
(closed) points of a variety $X$.  We say $\mc{C}$ is \emph{critically dense} in $X$ 
if every proper Zariski-closed subset $Y \subsetneq X$ contains only finitely many 
points of $\mc{C}$.  
\end{definition}

We may now prove our main result characterizing the noetherian property for $R$.  

\begin{theorem}
\label{R noetherian} Let $R = R(\varphi,c)$ for some $(\varphi,c) \in (\aut \mb{P}^t) 
\times \mb{P}^t$ such that Hypothesis~\ref{distinct} holds.  As always, set $c_i = 
\varphi^{-i}(c)$.  Then 
\begin{enumerate} 
\item \label{R noetherian 1} 
 $R(\varphi,c)$ is left noetherian if and only the set $\{c_i\}_{i \geq 0} $ is 
critically dense in $\mb{P}^t$. 
\item \label{R noetherian 2}  
$R(\varphi,c)$ is right noetherian if and only the set $\{c_i \}_{i \leq -1}$ is 
critically dense in $\mb{P}^t$. 
\item \label{R noetherian 3} 
$R(\varphi,c)$ is noetherian if and only the set $\{c_i\}_{i \in \mb{Z}}$ is 
critically dense in $\mb{P}^t$.
\end{enumerate}
\end{theorem}
\begin{proof} (1) Set $\mc{C} = \{c_i\}_{i \geq 0} $ and suppose that $\mc{C}$ is critically
dense.  Then for any nonzero homogeneous polynomial $f \in R$, the set $D = \{i \in 
\mb{N} | f(c_i) = 0 \}$ has finite cardinality, so by Propositions~\ref{1st type 
noetherian} and \ref{2nd type noetherian} the left $R$-modules $(Sf \cap R)/Rf$ and 
$S/(R + Sf)$ are noetherian.  By Lemma~\ref{filtration by types}, for any finitely 
generated homogeneous left ideal $I$ of $R$, the left $R$-module $(SI \cap R)/I$ is 
noetherian. By Lemma~\ref{noetherian passes down}, $R$ is a left noetherian ring.

Conversely, if $\mc{C}$ fails to be critically dense, then we may choose a nonzero 
homogeneous polynomial $h \in S$ which vanishes at infinitely many points of 
$\mc{C}$. Since by Lemma~\ref{essential} we know that $R \hookrightarrow S$ is an 
essential extension of $R$-modules, there exists a homogeneous $g \in R$ such that $0 
\not = f = gh \in R$. Then $f$ also vanishes at infinitely many points of $\mc{C}$, 
and so by Proposition~\ref{1st type noetherian}, the left $R$-module $(Sf \cap R)/Rf$ 
is not noetherian.  Since this module is a subfactor of $R$, we conclude that $R$ is 
not a left noetherian ring.

(2) Using Lemma~\ref{switch sides}(1), this part follows immediately from part (1).

(3) This follows from  the fact that for any infinite sets $\mc{C}_1, \mc{C}_2 \subseteq
\mb{P}^t$, $\mc{C}_1 \cup \mc{C}_2$ is critically dense if and only if both $\mc{C}_1$
and $\mc{C}_2$ are.
\end{proof}

In Section~\ref{examples} we will examine the 
critical density condition appearing in Theorem~\ref{R noetherian} more closely.  
In particular, we shall prove that there exist many choices of 
$\varphi$ and $c$ for which $R(\varphi,c)$ is noetherian:
\begin{proposition}
\label{diag case 1} (See Theorem~\ref{diag case} below) Let $\varphi$ be the 
automorphism of $\mb{P}^t$ defined by $(a_0: a_1: \dots : a_t ) \mapsto (a_0: p_1 
a_1: p_2 a_2: \dots : p_t a_t)$, and let $c$ be the point $(1:1:\dots:1) \in 
\mb{P}^t$.  If the scalars $\{p_1, p_2, \dots p_t \}$ are algebraically independent 
over the prime subfield of $k$, then $\{\varphi^i(c)\}_{i \in \mb{Z}}$ is critically 
dense and $R(\varphi,c)$ is noetherian. \hfill  $\Box$
\end{proposition}

The noetherian case is our main interest, so in the remainder 
of the paper (except \S\ref{examples})
we will assume the following hypothesis.
\begin{hypothesis} 
\label{noeth hyp} \label{cdc-index} 
Let $c_i = \varphi^{-i}(c)$.  Assume that $\varphi$ and $c$ are chosen such that the 
point set $\{c_i\}_{i \in \mb{Z}}$ is critically dense in $\mb{P}^t$, so that 
$R(\varphi,c)$ is noetherian.  We will refer to this as the \emph{critical density 
condition}.    
\end{hypothesis}
  
Below, we will frequently use the following exact sequence to study an arbitrary 
cyclic left $R$-module $R/I$: 
\begin{equation}
\label{exact seq}
 0 \ra (SI \cap R)/I \ra R/I \ra S/SI \ra S/(R + SI) \ra 0.
\end{equation}
We note for future reference what the results of this section tell us about the terms 
of this sequence.  
\begin{lemma}
\label{comparing cyclic R and S modules} Assume the critical density condition, and 
let $0 \neq I$ be a graded left ideal of $R$.  
\begin{enumerate} 
\item \label{comparing cyclic R and S modules 1} 
As left $R$-modules, $(SI \cap R)/I$ and  $S/(R + SI)$  
have finite filtrations with factors which are either torsion or a tail of the 
shifted $R$-point module $_R (\pointmod{c_{-1}})[-i]$ for some $i \geq 0$.  In 
particular, $S/(R + SI)$ is a noetherian left $R$-module. 
\item \label{comparing cyclic R and S modules 2} 
$_R (S/J)$ is a noetherian module for all nonzero left 
ideals $J$ of $S$.
\end{enumerate}
\end{lemma}
\begin{proof} 
(1) Let $0 \neq f \in R$ be arbitrary.  Since $\{c_i\}_{i \in \mb{Z}}$ is a 
critically dense set, $f(c_i) = 0$ holds for only finitely many $i \in \mb{Z}$.  Then 
by the results \ref{one point}---\ref{2nd type noetherian}, the left $R$-modules $(Sf 
\cap R)/Rf$ and $S/(R + Sf)$ are isomorphic to finite direct sums of shifted point 
modules of the form $_R (\pointmod{c_{-1}})[-i]$ for various $i \geq 0$. Now using 
Lemma~\ref{filtration by types}, it is clear that $(SI \cap R)/I$ has a filtration of 
the right kind.  Similarly, $S/(R + SI)$ is a homomorphic image of $S/(R + Sf)$ for 
any $0 \neq f \in I$, so it also has the required filtration and is clearly 
noetherian.  

(2) It is immediate from the exact sequence \eqref{exact seq} for $I = Rr$ and part 
(1) that $_R (S/Sr)$ is noetherian for any homogeneous $0 \neq r \in R$.  It is 
enough to show that $_R (S/Sx)$ is noetherian for an arbitrary homogeneous $0 \neq x 
\in S$.  There is some nonzero homogeneous $y \in R$ such that $yx \in R$, by 
Lemma~\ref{essential}.  Then since $(S/Syx)$ is a noetherian $R$-module, so is 
$S/Sx$. 
\end{proof}

%_____________________________________________________________________________
\section{Point modules and the strong noetherian property}
\label{point modules}

Let $S = S(\varphi)$ and $R = R(\varphi,c)$ for $(\varphi,c)$ satisfying the critical 
density condition, so that $R$ is noetherian.  Recall the definition of the strong 
noetherian property:
\begin{definition}
\label{strongnoetherian} A $k$-algebra $A$ is called \emph{strongly (left) 
noetherian} if $A \otimes_k B$ is a left noetherian ring for all commutative 
noetherian $k$-algebras $B$. 
\end{definition}
Artin and Zhang showed that the point modules for a strongly noetherian algebra have 
a nice geometric structure.  The following is a special case of their theorem.
\begin{theorem} \cite[Corollary~E4.11, Corollary~E4.12]{AZ}
\label{pm cor intro}
 Let $A$ be a connected $\mb{N}$-graded strongly noetherian
algebra over an algebraically closed field $k$.
\begin{enumerate}
\item \label{pm cor intro 1}
The point modules over $A$ are naturally parameterized by a commutative 
projective scheme over $k$.
\item \label{pm cor intro 2} 
There is some $d \geq 0$ such that every point module $M$ for $A$ is uniquely 
determined by its truncation $M/M_{\geq d}$. \hfill $\Box$
\end{enumerate} 
\end{theorem}

Using an explicit presentation for the ring, Jordan \cite{Jordan} classified the 
point modules for the algebra $R$ in a special case.   In this section, we classify 
the point modules for the rings $R(\varphi,c)$ in general using a different method 
which does not rely on relations and get a similar result.  The classification will 
show that part (2) of Theorem~\ref{pm cor intro} fails for $R$, and so $R$ cannot be 
strongly noetherian.  

In this section we will make frequent use of the criterion for $R$-membership given 
in Theorem~\ref{char of R} without comment. Also, recall that $\regmult$ indicates 
multiplication in the polynomial ring $U$, and juxtaposition indicates multiplication 
in $S$ (or $R$).   
Some of the results in this section will rely on the following technical commutative 
lemma which is proved in the appendix.
\begin{lemma} 
\label{case not on a line} (Lemma~\ref{case not on a line app}) Let the points $d_1, 
d_2, \dots, d_n, d_{n+1} \in \mb{P}^t$ be distinct, and assume that the points $d_1, 
\dots, d_n$ do not all lie on a line.  Let $\m_i \subseteq U$ be the homogeneous 
ideal corresponding to $d_i$.  
\begin{enumerate}
\item \label{case not on a line 1}
$(\bigcap_{i=1}^n \m_i)_{n-1} \circ (\m_{n+1})_1 = (\bigcap_{i = 1}^{n+1} \m_i)_n$. 
\item \label{case not on a line 2}
$(\bigcap_{i=1}^n \m_i)_{n-1} \circ (\m_1)_1 = (\bigcap_{i = 2}^n \m_i \cap \m_1^2)_n$. 
\item \label{case not on a line 3}
$(\bigcap_{i=2}^n \m_i \cap \m_1^2)_n \circ (\m_{n+1})_1 = (\bigcap_{i = 2}^{n+1} 
\m_i \cap \m_1^2)_{n+1}$. 
\item \label{case not on a line 4}
Let $b_1, b_2 \in \mb{P}^t$, with corresponding ideals $\mf{n_1}, \mf{n_2}$, be such 
that $b_j \neq d_i$ for $j = 1, 2$ and $1 \leq i \leq n$.  Then $(\bigcap_{i=1}^n 
\m_i \cap \mf{n}_1)_n = (\bigcap_{i=1}^n \m_i \cap \mf{n}_2)_n$ implies $b_1 =b_2$.  
\hfill $\Box$
\end{enumerate} 
\end{lemma} 

We have already seen that the point modules over $S$ are easily classified up to 
isomorphism---they are simply the $\{\pointmod{d} \mid d \in \mb{P}^t \}$ (recall 
Notation~\ref{point mod}).  There is a close relationship between the point modules 
over the rings $S$ and $R$, as we begin to see in the next proposition.
\begin{proposition}
\label{1 critical}
\begin{enumerate}
\item \label{1 critical 1} 
Let $M$ be a point module over $R$.  Then $M_{\geq n} \cong {}_R (P_{\geq n})$ for 
some $S$-point module $P$ and some $n \geq 0$. 
\item \label{1 critical 2} 
If $_R \pointmod{d_1}_{\geq n} \cong {}_R \pointmod{d_2}_{\geq n}$ for $d_1, d_2 \in 
\mb{P}^t$ and some $n \geq 0$, then $d_1 = d_2$. 
\end{enumerate} 
\end{proposition}
\begin{proof}
(1) We have $M = R/I$ for a unique point ideal $I$ of $R$.  We will use 
the exact sequence \eqref{exact seq}; there are two cases.

Suppose first that $(SI \cap R)/I = 0$.  Then we have an injection $R/I \ra S/SI$.  
By Lemma~\ref{comparing cyclic R and S modules}(1) we know that $\GK_R (S/(R+SI)) 
\leq 1$, and clearly $\GK_R (R/I) = 1$, so that $\GK_R (S/SI) = 1$ since GK-dimension 
is exact for modules over the graded noetherian ring $R$.  By Lemma~\ref{comparing 
cyclic R and S modules}(2), $_R (S/SI)$ is finitely generated, and so $\GK_S (S/SI) = 
\GK_R (S/SI) = 1$ since for finitely generated modules the GK-dimension depends only 
on the Hilbert function.  Now may choose a filtration of $S/SI$ composed of cyclic 
critical $S$-modules, where the factors must be 
shifts of $S$-point modules and $_S k$ (Lemma~\ref{cyclic critical filter}).  
Since $M$ is a $R$-submodule of $S/SI$, this 
forces some tail of $M$ to agree with a tail of an $S$-point module.

Suppose instead that $N = (SI \cap R)/I \neq 0$.  Then $N$ is a nonzero submodule of 
the point module $M$, so it is equal to a tail of $M$.  By Lemma~\ref{comparing 
cyclic R and S modules}(1), some tail of $N$, and thus a tail of $M$, must be 
isomorphic as an $R$-module to a tail of some $\pointmod{c_{-1}}[-i] \cong 
\pointmod{c_{-1-i}}_{\geq i}$ (using also Lemma~\ref{truncating point modules}).

(2) By Lemma~\ref{truncating point modules}, we have for any $d \in 
\mb{P}^t$ that $\pointmod{d}_{\geq n} \cong \pointmod{\varphi^n(d)}[-n]$ as 
$S$-modules.  Thus we may reduce to the case where $n = 0$.

Since $\ann_S \pointmod{d_i}_0 = \m_{d_i}$, we must have $\m_{d_1} \cap R = \m_{d_2} 
\cap R$.  In degree $m$ this means 
\begin{equation}
\label{comparing point ideals} (\m_{c_0} \cap \dots \cap \m_{c_{m-1}} \cap 
\m_{d_1})_m = (\m_{c_0} \cap \dots \cap \m_{c_{m-1}} \cap \m_{d_2})_m.
\end{equation}
Suppose first that $d_1, d_2 \not \in \{c_i\}_{i \in \mb{N}}$.  Since the point set 
$\{c_i\}_{i \in \mb{Z}}$ is critically dense, it follows that for $m \gg 0$ the 
points $\{c_i\}_{i = 0}^{m-1}$ do not all lie on a line.  Then by Lemma~\ref{case not 
on a line}(4), the equation \eqref{comparing point ideals} for $m \gg 0$ implies that 
$d_1 = d_2$.  

Otherwise we may assume without loss of generality that $d_1 = c_j$ for some $j \geq 
0$ and that $d_2 \not \in \{c_i\}_{i = 0}^{j-1}$.  Then the equation \eqref{comparing 
point ideals} for $m = j+1$ violates Lemma~\ref{Hilbert of point ideal} unless $d_1 = 
c_j = d_2$. 
\end{proof}

We may now classify the point modules over the ring $R(\varphi,c)$.  
\begin{theorem}
\label{classifying point ideals} Assume the critical density condition 
(Hypothesis~\ref{noeth hyp}).   
\begin{enumerate} 
\item \label{classifying point ideals 1} 
For any point $d \in \mb{P}^t 
\setminus \{c_i\}_{i \geq 0}$, the $S$-point module $\pointmod{d}$ is an $R$-point 
module, with point ideal $(R \cap \m_d)$.  
\item \label{classifying point ideals 2} 
For each $i \geq 0$, the $S$-module $\pointmod{c_i}_{\geq i+1}$ is a shifted 
$R$-point module.  There is a $\mb{P}^{t-1}$-parameterized family of non-isomorphic 
$R$-point modules \label{Pcie-index} $\{ \pointmod{c_i, e} \mid e \in \mb{P}^{t-1} 
\}$ with $\pointmod{c_i, e}_{\geq i+1} \cong {}_R \pointmod{c_i}_{\geq i+1}$ and 
$\pointmod{c_i, e}_{\leq i} \cong{} _R \pointmod{c_i}_{\leq i}$ for any $e \in 
\mb{P}^{t-1}$.  These are exactly the point modules whose point ideals contain the 
left ideal $(R \cap \m_{c_i}^2)$ of $R$.  
\item \label{classifying point ideals 3} 
All of the point modules given in parts (1) and (2) above are non-isomorphic, and 
every point module over $R(\varphi,c)$ is isomorphic to one of these.
\end{enumerate}
\end{theorem}
\begin{proof}
Suppose that $d \in \mb{P}^t$, so $\pointmod{d} = S/\m_d$ by definition.  For any $i 
\geq 0$,
\begin{equation}
\label{when breaks off}
R_1 (\pointmod{d})_i = 0 \Llra R_1 S_i \subseteq \m_d \Llra (\m_{c_i})_1 \regmult U_i \subseteq \m_d \Llra d = c_i.  
\end{equation}  

(1) Let $d \not \in \{c_i\}_{i \geq 0}$.  
In this case it is clear from \eqref{when breaks off} that $\pointmod{d}$ is already an $R$-point
module.  Also, the corresponding point ideal is $\ann_R \pointmod{d}_0 = R \cap \m_d$. 

(2) Fix some $i \geq 0$.  
From \eqref{when breaks off} it is clear that $M ={} _R \pointmod{c_i} = M_{\leq i} \oplus M_{\geq i+1}$
where $M_{\leq i}$ is the torsion submodule of $M$ and $M_{\geq i+1}$ is a shifted $R$-point module.

We define a left ideal $J = J^{(i)}$ of $R$ by setting $J_{\leq 
i} = (R \cap \mf{m}_{c_i})_{\leq i}$ and $J_{\geq i+1} = (R \cap 
\mf{m}_{c_i}^2)_{\geq i+1}$.  To check that $J$ really is a left ideal of $R$, one calculates
\[
R_1 J_i = \phi^i(\m_{c_0})_1 \circ J_i = (\mf{m}_{c_i})_1 \circ (R_i \cap 
\mf{m}_{c_i})_i \subseteq R_{i+1} \cap \mf{m}_{c_i}^2 = J_{i+1}.
\]

We will now classify the point ideals of $R$ which contain $J$. By Lemma~\ref{Hilbert 
of point ideal}, the Hilbert function of $R/J$ must be  
\[ \dim_k (R/J)_n =
\begin{cases}
1 & n \leq i \\
t & n \geq i+1.
\end{cases}
\]
Then using Lemma~\ref{Hilbert of point ideal} again, the natural injection 
\[
(R/J)_{\geq i+1} = 
\frac{(\m_{c_0} \cap \m_{c_1} \cap \dots \cap \m_{c_i})_{\geq i+1}}
{(\m_{c_0} \cap \m_{c_1} \cap \dots \cap \m_{c_i}^2)_{\geq i +1}}  
\hookrightarrow (\mf{m}_{c_i} / \mf{m}_{c_i}^2)_{\geq i+1}
\] is an 
isomorphism of left $R$-modules, since the Hilbert functions on both sides are the 
same.

As a module over the polynomial ring $U$, we have an isomorphism 
\[
(\mf{m}_{c_i} /
\mf{m}_{c_i}^2)_{\geq i+1} \cong \bigoplus_{j=1}^{t} (U/\mf{m}_{c_i})_{\geq
i+1}
\]
which by the equivalence of categories $U\Gr \sim S\Gr$ 
translates to an $S$-isomorphism as follows:   
\[
_S (\mf{m}_{c_i} / \mf{m}_{c_i}^2)_{\geq i+1} \cong \bigoplus_{j=1}^{t}
\pointmod{c_i}_{\geq i+1}.
\]
By part (1), $\pointmod{c_i}_{\geq i+1} \cong \pointmod{c_{-1}}[-i-1]$ is a shifted 
$R$-point module, so we conclude that $M = (R/J)_{\geq i+1}$ is a direct sum of $t$ 
isomorphic shifted $R$-point modules.  Then every choice of a codimension-one vector 
subspace $V = L/(J_{i+1})$ of $(R/J)_{i+1}$ generates a different $R$-submodule $N$ 
of $M$ with $M/N \cong {}_R \pointmod{c_i}_{\geq i+1}$, and then $J + RL$ is a point 
ideal for $R$.  Clearly any point ideal containing $J$ must arise in this way, and 
the set of codimension-one subspaces of $(R/J)_{i+1}$ is parameterized by 
$\mb{P}^{t-1}$.  Thus the set of point ideals of $R$ which contain $J$ is naturally 
parameterized by a copy of $\mb{P}^{t-1}$.

For each $e \in \mb{P}^{t-1}$, we have a corresponding point ideal $I$ containing $J$ 
and we set $P(c_i ,e) = R/I$.  Then $P(c_i,e)_{\geq i+1} \cong{} _R 
\pointmod{c_i}_{\geq i+1}$ and $P(c_i,e)_{\leq i} \cong{} (R/J)_{\leq i} \cong {}_R 
\pointmod{c_i}_{\leq i}$.   

Finally, note that all of the point ideals constructed above contain $(R \cap 
\m_{c_i}^2)$.  Conversely, if $I$ is any point ideal which contains $(R \cap 
\m_{c_i}^2)$, then $I$ contains $J$, since $J/(R \cap \m_{c_i}^2)$ is torsion and 
$I$, being a point ideal, is closed under extensions inside $R$ by torsion modules.  
Then $I$ is one of the point ideals we already constructed.  
This finishes the proof of Part (2).

(3) Note that for fixed $i$ the $P(c_i, e)$ are non-isomorphic for distinct $e$ by 
construction; then it follows easily from Proposition~\ref{1 critical}(2) that all of the 
point modules we have constructed in parts (1) and (2) are non-isomorphic.

Suppose that $M$ is an $R$-point module.  Let $\mc{P}$ be the set of all $R$-modules 
isomorphic to a shift of one of the $R$-point modules constructed in parts (1) and 
(2) above.  By Proposition~\ref{1 critical}, $M_{\geq n} \cong{} _R 
\pointmod{d}_{\geq n}$ for some $n \geq 0$ and $d \in \mb{P}^t$.  For $m \gg 0$, note 
that $\varphi^m(d) \not \in \{c_i\}_{i \geq 0}$ and so $M_{\geq m+n} \in \mc{P}$ by 
part (1) above.  Thus to finish the proof of (3) it is enough by induction to show 
that given any $R$-point module $N$, if $N_{\geq 1} \in \mc{P}$ then $N \in \mc{P}$. 

Let $N$ be an $R$-point module such that $N_{\geq 1} \in \mc{P}$.  Let $I = \ann_R 
N_0$ be the point ideal of $N$.  There are a number of cases.  

\emph{Case 1}.  Suppose first that $N_{\geq 1} \cong{} _R \pointmod{d}[-1]$ for some 
$d \not \in \{c_i\}_{i \geq -1}$.  Then $(R \cap \m_d) R_1 \subseteq I$, in other 
words 
\[ 
(\m_{\varphi^{-1}(d)} \cap \m_{c_1} \cap \m_{c_2} \cap \dots \cap \m_{c_n})_n \circ 
(\m_{c_0})_1 \subseteq I_{n+1} 
\]
for each $n \geq 0$.  By the critical density condition, for $n \gg 0$ the points 
$\{c_1, \dots c_n\}$ will not lie on a line, so that Lemma~\ref{case not on a 
line}(1) applies and gives 
\[ 
(\m_{\varphi^{-1}(d)} \cap \m_{c_0} \cap \m_{c_1} \cap \m_{c_2} \cap \dots \cap 
\m_{c_n})_{n+1} \subseteq I_{n+1}.
\]
In other words, $(R \cap \m_{\varphi^{-1}(d)})_{\geq m} \subseteq I$ for $m \gg 0$. 
Note that $\varphi^{-1}(d) \not \in \{c_i\}_{i \geq 0}$, so that $(R \cap 
\m_{\varphi^{-1}(d)})$ is one of the point ideals appearing in part (1) above.  Since 
$I$ is a point ideal and is thus closed under extensions inside $R$ by torsion 
modules, $(R \cap \m_{\varphi^{-1}(d)}) \subseteq I$ and so comparing Hilbert 
functions, $(R \cap \m_{\varphi^{-1}(d)}) = I$.  Thus $N \cong{} _R 
\pointmod{\varphi^{-1}(d)} \in \mc{P}$.

\emph{Case 2}.  If $N_{\geq 1} \cong{} _R \pointmod{c_{-1}}$, then $(R \cap 
\m_{c_{-1}}) R_1 \subseteq I$; using Lemma~\ref{case not on a line}(2) and a similar 
argument to that in case 1, this implies that $(R \cap \m_{c_0}^2)_{\geq m} \subseteq 
I$ for $m \gg 0$.  Then since $I$ is a point ideal, $(R \cap \m_{c_0}^2) \subseteq 
I$. By part (2) above this forces $N \cong P(c_0, e)$ for some $e \in \mb{P}^{t-1}$, 
and so $N \in \mc{P}$.
   
\emph{Case 3}. Let $N_{\geq 1} \cong{} P(c_i, e)[-1]$ for some $i \geq 0$ and $e \in 
\mb{P}^{t-1}$.  Then $(R \cap \m_{c_i}^2) R_1 \subseteq I$.  Now 
the same argument as in the other cases, except using Lemma~\ref{case not on a line}(3), 
will show that $N \cong P(c_{i+1}, e') \in \mc{P}$ for some $e' \in \mb{P}^{t-1}$.
\end{proof}

The failure of the strong noetherian property for $R = R(\varphi,c)$ now follows 
immediately from Theorem~\ref{classifying point ideals}(2).  This proves 
Theorem~\ref{main result 1} from the introduction. 
\begin{theorem}
\label{not strongly noetherian} Assume the critical density condition.  Then $R = 
R(\varphi,c)$ is a connected graded noetherian algebra, finitely generated in degree 
$1$, which is noetherian but not strongly noetherian. 
\end{theorem}
\begin{proof}
We only need to prove that $R$ is not strongly noetherian.  For each $i \geq 0$, 
Theorem~\ref{classifying point ideals}(2) provides a whole $\mb{P}^{t-1}$ of point 
modules $P(c_i, e)$ which have isomorphic truncations $P(c_i, e)_{\leq i}$.  By 
Theorem~\ref{pm cor intro}(2), $R$ cannot be a strongly noetherian $k$-algebra.
\end{proof}

We remark that the point modules over $R$ still appear to have an interesting 
geometric structure.  By Theorem~\ref{classifying point ideals}, there is a single 
point module corresponding to each point $d \in \mb{P}^t \setminus \{c_i\}_{i \geq 
0}$ and a $\mb{P}^{t-1}$-parameterized family of exceptional point modules 
corresponding to each point $c_i$ with $i \geq 0$.  Since blowing up $\mb{P}^t$ at a 
point in some sense  replaces that point by a copy of $\mb{P}^{t-1}$, the intuitive 
picture of the geometry of the point modules for $R$ is an infinite blowup of 
projective space at a countable point set.  

%________________________________________________
\section{Extending the base ring}
\label{explicit ring} 
Let $S = S(\varphi)$ and $R = R(\varphi,c)$ for $(\varphi,c)$ 
satisfying the critical density condition, and let $c_i = \varphi^{-i}(c) \in 
\mb{P}^t_k$ as usual.  We now know by Theorem~\ref{not strongly noetherian} that $R$ 
is not strongly noetherian, but this proof is quite indirect and it is not obvious 
which choice of extension ring $B$ makes $R \otimes_k B$ non-noetherian.  In this 
section we construct such a noetherian commutative $k$-algebra $B$ which is even a 
UFD.

Let $B$ be an arbitrary commutative $k$-algebra which is a domain. We will use 
subscripts to indicate extension of the base ring, so for example we write 
$R_B= R \otimes_k B$. The automorphism $\phi$ of $U$ 
extends uniquely to an automorphism of $U_B$ fixing $B$, which we also call $\phi$.  
We continue to identify the underlying $B$-module of $S_B$ with that of $U_B$, and we 
use juxtaposition for multiplication in $S_B$ and the symbol $\circ$ for 
multiplication in $U_B$, as in our current convention (see \S\ref{Zhang twists}).  
The multiplication of $S_B$ is still given by $fg = \phi^n(f) \regmult 
g$ for $f \in (S_B)_m$, $g \in (S_B)_n$; in other words, $S_B$ is the left Zhang 
twist of $U_B$ by the twisting system $\{\phi^i\}_{i \in \mb{N}}$, just as before.  

Let $d$ be a point in $\mb{P}^t_k$.  Since the homogeneous coordinates for $d$ are 
defined only up to a scalar multiple in $k^{\times}$, given $f \in U_B$ the expression
$f(d)$ is defined up to a nonzero element of $k$; we will use this notation only in 
contexts where the ambiguity does not matter.  For example, the condition $f(d) = 0$ 
makes sense and is equivalent to the condition $f \in \m_d \circ U_B$, where $\m_d 
\subseteq U$ and $\m_d \circ U_B$ is a graded prime ideal of $U_B$.  

The natural analog of Theorem~\ref{char of R} still holds in this setting: 
\begin{proposition}
\label{char of R_B}
For all $n \geq 0$, we have 
\[ (R_B)_n = \{f \in (U_B)_n\ \text{such that}\ f(c_i) = 0\ \text{for}\ 0 \leq i \leq n-1 \}. 
\]
\end{proposition}
\begin{proof}
As subsets of $U_B$, using Theorem~\ref{char of R} we have 
\[
(R_B)_n = R_n \otimes B = (\cap_{i = 0}^{n-1} \m_{c_i})_n \otimes B = \cap_{i = 
0}^{n-1} (\m_{c_i} \circ U_B)_n 
\] 
and the proposition follows.
\end{proof}

We now give sufficient conditions on $B$ for the ring $R \otimes_k B$ to fail to have the left noetherian property.
\begin{proposition}
\label{when not noeth} Assume that $B$ is a UFD. Suppose that there exist nonzero 
homogeneous elements $f, g \in (U_B)_1$ satisfying the following conditions: 
\begin{enumerate}
\item \label{when not noeth 1}
$f(c_i)$ divides $g(c_i)$ in $B$ for all $i \geq 0$. 
\item \label{when not noeth 2}
For all $i \gg 0$, $f(c_i)$ is not a unit of $B$. 
\item \label{when not noeth 3}
$\gcd(f, g) = 1$ in $U_B$. 
\end{enumerate} 
Then $R \otimes_k B$ is not a left noetherian ring.
\end{proposition}
\begin{proof}
Note that $U_B \cong B[x_0, x_1, \dots, x_t]$ is a UFD, since $B$ is, so condition 
(3) makes sense.

For convenience, fix some homogeneous coordinates for the $c_i$. For each $n \geq 0$, 
we may choose a polynomial $\theta_n \in S_n$ with coefficients in $k$ such that 
$\theta_n(c_i) = 0$ for $-1 \leq i \leq n-2$ and $\theta_n(c_{n-1}) \neq 0$.  This is 
possible, for example, by Lemma~\ref{Hilbert of point ideal}. By hypothesis~(1), for 
each $n \geq 0$ we may write $\Omega_n = g(c_n)/f(c_n) \in B$. Now let $t_n = 
\theta_n(\Omega_n f - g) \in (S_B)_{n+1}$ for each $n \geq 0$.

Since $\phi(\theta_n)$ vanishes at $c_i$ for $0 \leq i \leq n-1$ and $[\Omega_n f 
-g](c_n) = 0$, the element $t_n = \phi(\theta_n) \circ (\Omega_n f - g)$ is in 
$(R_B)_{n+1}$, by Proposition~\ref{char of R_B}. We will show that for $n \gg 0$ we 
have $t_{n+1} \not \in \sum_{i = 0}^n (R_B) t_i$, which will imply that $R_B$ is not 
left noetherian.

Suppose that $t_{n+1} = \sum_{i = 0}^n r_i t_i$ for some $r_i \in (R_B)_{n+1-i}$.  
Writing out the explicit expressions for the $t_i$, this is 
\[
\theta_{n+1} (\Omega_{n+1} f -g) = \sum_{i=0}^n r_i \theta_i (\Omega_i f -g).
\]
Considering these expressions in $U_B$, after some rearrangement we obtain (since 
$f,g$ have degree 1)
\[
\phi[\theta_{n+1} \Omega_{n+1} - \sum_{i=0}^n r_i \theta_i \Omega_i] \circ f +
\phi[-\theta_{n+1} + \sum_{i =0}^n r_i \theta_i ] \circ g = 0.
\]
Now by hypothesis~(3), $g$ must divide the polynomial 
\[
h = \phi[\theta_{n+1} \Omega_{n+1} - \sum_{i=0}^n r_i \theta_i \Omega_i] = 
\Omega_{n+1} \phi(\theta_{n+1})  - \sum_{i=0}^n \Omega_i \phi^{i+1}(r_i) \circ 
\phi(\theta_i).
\]
We note that $[\phi(\theta_{n+1})](c_{n+1}) \in k^{\times}$ by the definition 
of the $\theta_i$, and $[\phi^{i+1}(r_i)](c_{n+1}) = 0$ for 
$0 \leq i \leq n$, since $r_i \in (R_B)_{n-i+1}$.  Thus evaluating at $c_{n+1}$ we 
conclude that $h(c_{n+1}) \in \Omega_{n+1}k^{\times}$ and thus 
$g(c_{n+1})$ divides $\Omega_{n+1}$.  But since $\Omega_{n+1} = g(c_{n+1})/f(c_{n+1})$, this 
implies that $f(c_{n+1})$ is a unit in $B$. For all $n \gg 0$, this contradicts 
hypothesis~(2), and so  $t_{n+1} \not \in \sum_{i = 0}^n (R_B) t_i$ for $n \gg 0$, as 
we wished to show.
\end{proof}

Next, we construct a commutative noetherian ring $B$ which satisfies the hypotheses
of Proposition~\ref{when not noeth}.  
We shall obtain such a ring as an infinite blowup of 
affine space, to be defined presently.  See \cite[Section 1]{ASZ} for more details 
about this construction.

Let $A$ be a commutative domain, and let $X$ be the affine scheme $\spec A$.  Suppose
that $d$ is a closed nonsingular point of $X$ with corresponding maximal ideal 
$\mf{p} \subseteq A$, and let $z_0, z_1, \dots, z_r$ be some choice of generators of 
the ideal $\mf{p}$ such that $z_0 \not \in \mf{p}^2$.  The \label{affineblowup-index} 
\emph{affine blowup} of $X$ at 
$d$ (with denominator $z_0$) is $X' = \spec A'$ where $A' = A[z_1 z_0^{-1}, z_2 
z_0^{-1}, \dots, z_r z_0^{-1}]$. 

Consider the special case where $A = k[y_1, y_2, \dots, y_t]$ is a polynomial ring, 
$X = \mb{A}^t$, and $d = (a_1, a_2, \dots, a_t)$.  The affine blowup of $\mb{A}^t$ at 
$d$ with the denominator $(y_1 - a_1)$ is $X' = \spec A'$ for the ring 
\[
A' = A[(y_2 - a_2)(y_1 -a_1)^{-1}, \dots, (y_t - a_t)(y_1 - a_1)^{-1}].
\]
Note that also $A' = k[y_1, (y_2 - a_2)(y_1 -a_1)^{-1}, \dots, (y_t - a_t)(y_1 - 
a_1)^{-1}]$, so $A'$ is itself isomorphic to a polynomial ring in $t$ variables over 
$k$ and $X' = \mb{A}^t$ as well.  The blowup map $X' \ra X$ is an isomorphism outside 
of the closed set $\{ y_1 = a_1 \}$ of $X$.

Given a sequence of points $\{ d_i = (a_{i1}, a_{i2}, \dots, a_{it}) \}_{i \geq 0}$ 
such that $a_{i1} \neq a_{j1}$ for $i \neq j$, we may iterate the blowup 
construction, producing a union of commutative domains 
\[
A \subseteq A_0 \subseteq A_1 \subseteq A_2 \subseteq \dots 
\]
where each $A_i$ is isomorphic to a polynomial ring in $t$ variables over $k$.  In 
this case we set $B = \bigcup A_i$ and $Y = \spec B$, and call $Y$ (or $B$) the 
\emph{infinite blowup} \label{infiniteblowup-index} of $\mb{A}^t$ at the sequence of points $\{d_i\}$.  
Explicitly, $B = A[\{(y_j - a_{ij})(y_1 - a_{i1})^{-1} | 2 \leq j \leq t,\ i \geq 0 
\} ]$. 

That there should be some connection between such infinite blowups and the algebras 
$R(\varphi,c)$ is strongly suggested by the following result (compare Theorem~\ref{R 
noetherian}).
\begin{theorem}
\label{inf blowup noeth} \cite[Theorem 1.5]{ASZ} The infinite blowup $B$ is a 
noetherian ring if and only if the set of points $\{d_i\}_{\geq 0}$ is a critically 
dense subset of $\mb{A}^t$. \hfill $\Box$
\end{theorem}

Now we show the failure of the strong noetherian property for the noetherian rings 
$R(\varphi,c)$ explicitly.  Below, we we will identify automorphisms of $\mb{P}^t$ 
with elements of $\pgl_{t+1}(k) = \gl_{t+1}(k) /k^{\times}$ \cite[page 151]{H}; 
explicitly, we let matrices in $\gl_{t+1}(k)$ act on the left on the homogeneous 
coordinates $(a_0:a_1: \dots :a_t)$ of $\mb{P}^t$, considered as column vectors. 

\begin{theorem}
\label{explicit B} Let $(\varphi,c) \in (\aut \mb{P}^t_k) \times \mb{P}^t_k$ satisfy 
the critical density condition.  There is an affine patch $\mb{A}^t \subseteq 
\mb{P}^t$ such that $\{c_i\}_{i \in \mb{Z}} \subseteq \mb{A}^t$.  Let $B$ be the 
infinite blowup of $\mb{A}^t$ at the points $\{c_i\}_{i \geq 0}$.  Then $R = 
R(\varphi, c)$ is noetherian, but $B$ is a commutative noetherian $k$-algebra which 
is a UFD such that $R \otimes_k B$ is not a left noetherian ring.
\end{theorem}
\begin{proof}

By changing coordinates, we may replace $\varphi$ by a conjugate without loss of 
generality, so we may assume that when represented as a matrix $\varphi$ is lower 
triangular.  Also, we may multiply this matrix by a nonzero scalar without changing 
the automorphism of $\mb{P}^t$ it represents, and so we also assume that the top left 
entry of the matrix is $1$.

By assumption the set of points $\{c_i\}_{i \in \mb{Z}}$ is critically dense in 
$\mb{P}^t$, and so $R(\varphi,c)$ is noetherian. Let $X_0$ be the hyperplane $\{x_0 = 
0\}$ of $\mb{P}^t$.  Since $\varphi$ is upper triangular, $\varphi(X_0) = X_0$, so if 
some $c_i \in X_0$ then $\{c_i\}_{i \in \mb{Z}} \subseteq X_0$ which contradicts the 
critical density condition.  So certainly $\{c_i\}_{i \in \mb{Z}} \subseteq \mb{A}^t 
= \mb{P}^t \setminus X_0$.  Since the top left entry of $\varphi$ is $1$, we may fix 
homogeneous coordinates for the $c_i$ of the form $c_i = (1: a_{i1}: a_{i2}: \dots : 
a_{it})$.  Let $y_i = x_i/x_0$, so that $k[y_1, y_2, \dots y_t]$ is the coordinate 
ring of $\mb{A}^t$.  In affine coordinates, $c_i = (a_{i1}, a_{i2}, \dots, a_{it})$.  

If $a_{i1} = a_{j1}$ for some $i < j$, then since $\varphi$ is lower triangular it 
follows that $a_{i1} = a_{k1}$ for all $k \in (j-i)\mb{Z}$.  Then the hyperplane 
$\{a_{i1} x_0 - x_1 = 0\}$ of $\mb{P}^t$ contains infinitely many of the $c_i$, again 
contradicting the critical density of $\{c_i\}_{i \in \mb{Z}}$.  So the scalars 
$\{a_{i1}\}_{i \in \mb{Z}}$ are all distinct, and the infinite blowup $B$ of 
$\mb{A}^t$ at the points $\{c_i\}_{i \geq 0}$ is well defined.  The ring $B$ is 
generated over $k[y_1, y_2, \dots, y_t]$ by the elements $\{ (y_j - a_{ij})(y_1 - 
a_{i1})^{-1} \mid 2 \leq j \leq t,\ i \geq 0 \}$.  Clearly the points $\{c_i\}_{i 
\geq 0}$ must be critically dense subset of $\mb{A}^t$, since they are a critically dense
subset of $\mb{P}^t$.  Thus $B$ is 
noetherian by Theorem~\ref{inf blowup noeth}. 

The ring $B$ is obtained as a directed union of $k$-algebras $A_i$ which are each isomorphic 
to a polynomial ring.  In each ring $A_i$ the group of units is just $k^{\times}$, 
and so this is also the group of units of $B$.  It follows that if $z \in A_i$ is an 
irreducible element of $B$, then $z$ is irreducible in $A_i$.  Since $B$ is 
noetherian, every element of $B$ is a finite product of irreducibles, and the 
uniqueness of such a decomposition follows by the uniqueness in each UFD $A_i$.  Thus 
$B$ is a UFD.  

Fix the two elements $f = y_1 x_0 -  x_1$ and $g = y_2 x_0 - x_2$ of $U_B \cong 
B[x_0, x_1, \dots, x_t]$.  Since $f$ and $g$ are homogeneous of degree $1$ in the 
$x_i$ and are not divisible by any non-unit of $B$, it is clear that $f$ and $g$ are 
distinct irreducible elements of $U_B$, and so in particular $\gcd(f, g) = 1$.  Now 
$f(c_i) = y_1 - a_{i1}$ and $g(c_i) = y_2 - a_{i2}$, so $\Omega_i = g(c_i)/f(c_i) = 
(y_2 - a_{i2})(y_1 - a_{i1})^{-1} \in B$ and thus $f(c_i)$ divides $g(c_i)$ for all $i \geq 
0$.

Finally, $f(c_i) = (y_1 - a_{i1})$ is not in the group of units $k^{\times}$ of $B$.  
We see that all of the hypotheses of Proposition~\ref{when not noeth} are satisfied, 
and so $R \otimes_k B$ is not left noetherian.
\end{proof}

%_______________________________________________________________________________________________
\section{Special subcategories and homological lemmas}
\label{homological lemmas} 

Let $S = S(\varphi)$ and $R = R(\varphi,c)$, and assume the critical density 
condition (Hypothesis~\ref{noeth hyp}), in particular that $R$ is noetherian. First, 
we introduce some notation for the subcategories of $S\Gr$ and $R\Gr$ which are 
generated by the ``distinguished'' $S$-point modules $P(c_i)$.  
\begin{definition}
\label{Sdist-index} \label{Rdist-index} \label{spec cats} 
\begin{enumerate} 
\item \label{spec cats 1} 
Let $S\cat$ be the full subcategory of $S\gr$ consisting of 
all $S$-modules $M$ with a finite $S$-module filtration whose factors are either 
torsion or a shift of $\pointmod{c_i}$ for some $i \in \mb{Z}$. 
\item \label{spec cats 2}
Let $R\cat$ be the full subcategory of $R\gr$ consisting of all $R$-modules 
$M$ having a finite $R$-module filtration whose factors are either torsion or a shift 
of the module $_R \pointmod{c_i}$ for some $i \in \mb{Z}$. 
\end{enumerate}
\end{definition}

Note that by Theorem~\ref{classifying point ideals}(1),(2), $_R \pointmod{c_i}$ is 
finitely generated for any $i \in \mb{Z}$ and so part (2) of the definition makes 
sense.  We also define $S\Cat$ to be the smallest full subcategory of $S\Gr$ 
containing $S\cat$ and closed under direct limits.  The subcategory $R\Cat$ of $R\Gr$ 
is defined similarly.  We will use frequently later in this section 
the fact that $S/R \in R\Cat$, which follows from Corollary~\ref{structure of S/R}.  

If $\mc{C}$ is any abelian category, a full subcategory $\mc{D}$ of $\mc{C}$ is 
called \emph{Serre} if for any short exact sequence \mbox{$0 \ra M' \ra M \ra M'' \ra 
0$} in $\mc{C}$, $M \in \mc{D}$ if and only if both $M' \in \mc{D}$ and $M'' \in 
\mc{D}$.  All of the subcategories defined above are clearly Serre.  In fact, we may 
describe $S\cat$ as the smallest Serre subcategory of $S\gr$ which contains all of 
the $P(c_i)[j]$; a similar description holds for the other categories.

The special categories over the two rings are related as follows. 
\begin{lemma}
\label{relation between special categories} Let $M \in S\gr$.  Then $_R M \in R\cat$ 
if and only if $_S M \in S\cat$.  Similarly, if $M \in S\Gr$ then $_R M \in R\Cat$ if 
and only if $_S M \in S\Cat$.
\end{lemma}
\begin{proof}
If $_S M \in S\cat$, then it follows directly from Definition~\ref{spec cats} that 
$_R M \in R\cat$.   Conversely, suppose that $_R M \in R\cat$.  Clearly $\GK_R (M) 
\leq 1$, so we have $\GK_S (M) \leq 1$ since we can measure GK-dimension using the 
Hilbert function.  By Lemma~\ref{cyclic critical filter}, $M$ has a finite filtration 
over $S$ with cyclic critical factors, which must in this case be shifts of $_S k$ 
and $S$-point modules.  Suppose that a shift of $\pointmod{d}$ is one of the factors 
occurring.  Then $N = {}_R \pointmod{d} \in R\cat$.  By the definition of $R\cat$, 
some tail of $N$ is isomorphic to a shift of some $_R \pointmod{c_i}$ for some $i \in 
\mb{Z}$.  Using Lemma~\ref{truncating point modules}, this forces $_R \pointmod{d} 
\cong {}_R \pointmod{c_j}$ for some $j \in \mb{Z}$ and so by Proposition~\ref{1 
critical}(2) we have $d = c_j$.  Thus the only point modules which may occur as 
factors in the $S$-filtration of $M$ are shifts of the $\pointmod{c_j}$ for $j \in 
\mb{Z}$ and so $M \in S\cat$.  

The second statement is an easy consequence of the definitions of $S\Cat$ and $R\Cat$ 
and the first statement.
\end{proof}

In the rest of this section, we gather some definitions and lemmas concerning homological
algebra over the rings $R$ and $S$.  
Let $A$ be a connected $\mb{N}$-graded $k$-algebra, finitely generated in degree $1$, 
and let $k = (A/A_{\geq 1})$.  We say that $A$ satisfies $\chi_i$ if $\dim_k 
\uext^j(k, M)< \infty$ for all $M \in A\gr$ and all $0 \leq j \leq i$, and that $A$ 
satisfies $\chi$ if $A$ satisfies $\chi_i$ for all $i \geq 0$.  If $M \in A\gr$, the 
\emph{grade} \label{grade-index} of $M$ is the number $j(M) = \min\{ i | 
\uext^i_A(M,A) \neq 0 \}$.  We say that $A$ is \emph{Cohen-Macaulay} \label{CM-index} 
if $j(M) + \GK(M) = \GK(A)$ for all $M \in A\gr$.   Finally, $A$ is 
\emph{Artin-Schelter regular} (or \emph{AS-regular}) if $A$ has finite global 
dimension $d$, finite GK-dimension, and satisfies the \label{Gorenstein-index} 
\emph{Gorenstein condition}: $\uext^i_A(_A k, A)= 0$ if $i \neq d$, and $\uext^d_A(_A 
k, A) \cong k_A{}$ (up to some shift of grading). 

The ring $S$ obtains many nice homological properties simply because it is a Zhang twist
of a commutative polynomial ring.  
\begin{lemma}
\label{CM for S}
\begin{enumerate}
\item \label{CM for S 1}
$S$ has global dimension $t+1$. 
\item \label{CM for S 2}
$S$ is Cohen-Macaulay. 
\item \label{CM for S 3}
$S$ is Artin-Schelter regular. 
\item \label{CM for S 4} 
$S$ satisfies $\chi$. 
\end{enumerate}
\end{lemma}
\begin{proof}
All of these properties are standard for the polynomial ring $U$.  Properties (1)-(3) 
follow for the Zhang twist $S$ of $U$ by \cite[Propositions~5.7, 5.11]{Zhang}.  Then 
since $S$ is Artin-Schelter regular it satisfies $\chi$ \cite[Theorem~8.1]{AZ94}. 
\end{proof}

Recall from \S\ref{Zhang twists} that there is an equivalence of categories $\theta: 
U\Gr \sim S\Gr$, and that we identify the graded left ideals of $S$ and the graded 
ideals of $U$.  Under the equivalence of categories we have $\theta(U/I) \cong 
S/I$.   Also, graded injective objects correspond under the equivalence and so it is 
immediate that $\ext^i_{U\Gr}(M, N) \cong \ext^i_{S\Gr}(\theta{M}, \theta{N})$ as vector spaces 
for all $M, N \in U\Gr$ and all $i \geq 0$. On the other hand, the relationship 
between $\uext$ groups over the two rings is more complicated, since the shift 
functors in $S$ and $U$ do not correspond under the equivalence of categories.  We 
work out this relationship in detail for cyclic modules, which is the only case we 
will need.
\begin{lemma}
\label{shifts and equiv} Let $M = U/I$ and $N = U/J$ for some graded ideals $I$ and 
$J$ of $U$.  For any $n \in \mb{Z}$ we have 
\begin{enumerate}
\item \label{shifts and equiv 1} 
$\theta((U/\phi^{-n}(J))[n]) \cong (S/J)[n]$. 
\item \label{shifts and equiv 2}  
$\uext^i_S(S/I,S/J)_n \cong \uext^i_U (U/I, U/\phi^{-n}(J))_n\
    \text{as $k$-spaces.}$
\end{enumerate}
\end{lemma}
\begin{proof}
(1) Obviously the property of being cyclic is preserved by the equivalence of 
categories, and so $\theta((U/\phi^n(J))[n])$ is a cyclic $S$-module generated in 
degree $-n$.  Thus we need only show that that the annihilator in $S$ of a generator 
is the left $S$-ideal $J$, but this is immediate from the definition of the 
equivalence $\theta$.

(2) By definition, $\uext^i_U (U/I, U/\phi^{-n}(J))_n = \ext^i_{U\Gr}(U/I, 
(U/\phi^{-n}(J))[n])$ and $\uext^i_S(S/I,S/J)_n = \ext^i_{S\Gr}(S/I,(S/J)[n])$.  We know 
that $\theta(U/I) \cong S/I$, and $\theta((U/\phi^{-n}(J))[n]) \cong (S/J)[n]$ by 
part (1).  We are done since the $\ext$ groups in the graded category correspond 
under the equivalence of categories.  
\end{proof}

The next proposition shows that the critical density of the set $\{c_i\}_{i \in 
\mb{Z}}$, besides characterizing the noetherian property for $R$, also has 
implications for the homological properties of the $S$-point modules $\pointmod{c_i}$.
The proof of the following commutative lemma may be found in the appendix. 
\begin{lemma} (Lemma \ref{comm ext lemma app})
\label{comm ext lemma} Let $I$ and $J$ be homogeneous ideals of $U$.  There is some 
$d \geq 0$ such that for all $n \in \mb{Z}$ for which $U/(I + \phi^n(J))$ is bounded, 
$\uext_U^i(U/I, U/\phi^n(J))$ has right bound $\leq d$.  \hfill $\Box$
\end{lemma}
\begin{proposition}
\label{exts over S} Assume the critical density condition, and let $N \in S\gr$.  
\begin{enumerate} 
\item \label{exts over S 1} 
$\dim_k \uext_S^p (\pointmod{c_i},N) < \infty$ for $0 \leq p \leq t-1$ and any $i \in 
\mb{Z}$. 
\item \label{exts over S 2} 
Let $M \in S\cat$.  Then  $\dim_k \uext_S^p (M,N) < \infty$ for $0 \leq p \leq t-1$.
\end{enumerate} 
\end{proposition}
\begin{proof}
(1) Since $N$ is finitely generated, it is easy to see that each graded piece of $E = 
\uext_S^p(\pointmod{c_i},N)$ is 
finite dimensional over $k$.  So it is enough to show that $E$ is bounded.  Note that 
$E$ is automatically left bounded since $N$ is \cite[Proposition~3.1.1(c)]{AZ94}. 
It remains to show that $E$ is right bounded.  Using a finite filtration of $N$ by 
cyclic modules, one reduces quickly to the case where $N$ is cyclic, say $N = S/I$.  

In case $I = 0$, $E = \uext^p_S(\pointmod{c_i}, S) = 0$ for $0 \leq p \leq t-1$ by 
the Cohen-Macaulay property of $S$ (Lemma~\ref{CM for S}(2)).

Now assume that $I \neq 0$.  By Lemma~\ref{shifts and equiv}(3), we have for each $n 
\geq 0$ the $k$-space isomorphism 
\[
\uext_S^p (S/\m_{c_i},S/I)_n \cong \uext_U^p (U/\m_{c_i},U/\phi^{-n}(I))_n.
\]
Now $\phi^{-n}(I) \subseteq \m_{c_i}$, or equivalently $I \subseteq \m_{c_{i + n}}$, 
can hold for at most finitely many $n$, since the points $\{c_i\}_{i \in \mb{Z}}$ are 
critically dense.  Thus for $n \gg 0$ we have $\phi^{-n}(I) \not \subseteq \m_{c_i}$, 
and the module  $U/(\phi^{-n}(I) + \m_{c_i})$ is bounded.  By Lemma~\ref{comm ext lemma}, 
there is some fixed $d \geq 0$ such that $\uext_U^p 
(U/\m_{c_i},U/\phi^{-n}(I))_n = 0$ as long as $n \geq d$. We conclude that $\uext_S^p 
(S/\m_{c_i},S/I)_n = 0$ for $n \gg 0$, as we wish.

(2) Since $M \in S\cat$, we may choose a finite filtration of $M$ with factors which 
are shifts of the point modules $\pointmod{c_i}$ or $_S k$.  Since $S$ satisfies 
$\chi$ by Lemma~\ref{CM for S}(4), $\dim_k \uext_S^p(k, N) < \infty$ for all $p \geq 
0$, and now the statement follows by part (1).
\end{proof}

For an $\mb{N}$-graded algebra $A$, if $L$ and $N$ are $\mb{Z}$-graded right and left 
$A$-modules respectively, then the $k$-space $\tor_i^A(L,N)$ has a natural 
$\mb{Z}$-grading which we emphasize by using the notation \label{utor-index} 
$\utor_i^A(L,N)$.  To study homological algebra over $R$, we will generally try to 
reduce to calculations over the ring $S$.  In particular, we will often use the 
following convergent spectral sequence, which is valid for any graded modules $_R M$ 
and $_S N$ \cite[Equation~(2.2)]{StZh2}:
\begin{equation}
\label{spectral sequence} \uext^p_S(\utor^R_q(S, M), N) \underset{p}{\Rightarrow} 
\uext^{p+q}_R(M,N). 
\end{equation}
We also note for reference the 5-term exact sequence arising from this spectral 
sequence \cite[11.2]{Rotman}:
\begin{multline}
\label{5 term exact} 0 \ra \uext_S^1(S \otimes_R M, N) \ra \uext^1_R(M,N) \ra 
\uhom_S(\utor_1^R(S,M), N) \\
\ra \uext^2_S(S \otimes_R M, N) \ra \uext^2_R(M,N). 
\end{multline}

In order to make effective 
use of the spectral sequence, we need some information about $\utor$.  
\begin{lemma}
\label{tor from S} Fix some $M \in R\gr$.  Also, let $Q$ be the right point module of 
$R$ such that $(S/R)_R \cong \oplus_{i =1}^{\infty} Q[-i]$ (Corollary~\ref{structure 
of S/R}(\ref{structure of S/R 2}), applied to the right side, which is valid by 
Lemma~\ref{switch sides}(1)).  Then
\begin{enumerate}   
\item \label{tor from S 1} 
$\utor_q^R(S,M) \in S\cat$ for any $q \geq 1$.  If $M$ is torsion, then $\utor_q^R(S, 
M) \in S\cat$ for $q \geq 0$. 
\item \label{tor from S 2} 
$\dim_k \utor_q^R(Q, M) < \infty$ for $q \geq 1$.
\end{enumerate}
\end{lemma}
\begin{proof} 
(1) From the long exact sequence in $\utor^R_q( -, M)$ associated to the short exact 
sequence of $R$-bimodules $0 \ra R \ra S \ra S/R \ra 0$, we see that 
\begin{equation}
\label{tor fact 1} \utor^R_q(S, M) \cong \utor^R_q(S/R, M)
\end{equation}
as left $R$-modules, for all $q \geq 2$.  Also, at the bottom of the long exact 
sequence we have 
\begin{equation}
\label{tor fact 2} 
0 \to \utor_1^R(S, M) \to \utor_1^R(S/R, M) \to M \to \dots 
\end{equation}
Thus there is at least an injection of left $R$-modules $\utor^R_q(S,M) \ra 
\utor^R_q(S/R,M)$ for all $q \geq 1$.  Now computing $N = \utor_q^R(S/R, M)$ using a 
free resolution of $M$, it is a subfactor of some direct sum of copies of $(S/R) \in R\Cat$, so 
$N \in R\Cat$.  Then $\utor^R_q(S,M) \in R\Cat$ and thus in $S\Cat$ for $q \geq 1$, 
using Lemma~\ref{relation between special categories}.  But since we may calculate 
$\utor^R_q(S,M)$ using a resolution of $M$ by free modules of finite rank, we have 
$\utor^R_q(S,M) \in S\cat$ for $q \geq 1$.  

If $M$ is torsion, we need to show in addition that $\utor_0^R(S, M)$ is in $S\cat$. 
It is enough to show this for $M = k$, in which case one calculates immediately that 
$\utor_0^R(S, M) \cong S/SR_{\geq 1} \cong S/\m_{c_0} = P(c_0)$ which is obviously in 
$S\cat$.

(2) As in part (1), $N = \utor^R_q(S/R,M) \in R\Cat$ and $\utor_q^R(S, M)$ is in 
$S\cat$ and thus in $R\cat$, by Lemma~\ref{relation between special categories}, for 
all $q \geq 1$.  Since $M \in R\gr$, we get using \eqref{tor fact 1} and \eqref{tor 
fact 2} that $N \in R\gr$ and so $N \in R\cat$ for all $q \geq 1$.

Note that by the definition of $R\cat$, the Hilbert function of $N$ is forced to 
satisfy $\dim_k N_m < C$ for some constant $C$, all $m \geq 0$.  Then since $\utor$ 
commutes with direct sums \cite[Proposition~2.4(1)]{AZ94}, $N \cong \oplus_{i 
=1}^{\infty} \utor_q(Q,M)[-i]$ as graded vector spaces, and so we must have $\dim_k 
\utor_q^R(Q, M) < \infty$ for $q \geq 1$.  
\end{proof}

As an easy consequence of the spectral sequence, we may show that $R$ and $S$ have no 
nontrivial extensions by torsion modules in the category of $R$-modules. 
\begin{lemma} 
\label{ext k to R is 0} 
$\uext^1_R(_R k, R) = 0 = \uext^1_R(_R k, S)$.
\end{lemma}
\begin{proof} 
Consider the long exact sequence in $\uext_R(k,-)$ associated to the short exact 
sequence $0 \ra R \ra S \ra S/R \ra 0$:
\begin{equation}
\label{ex seq for ext k to R}  
\dots \ra \uhom_R(k, S/R) \ra \uext_R^1(k, R) \ra 
\uext_R^1(k, S) \ra \dots
\end{equation}
Now $_R(S/R)$ is torsionfree, since it is isomorphic to a direct sum of point modules 
by Corollary~\ref{structure of S/R}(2).  Thus $\uhom_R(k, S/R) = 0$. 

To analyze the group $\uext_R^1(k, S)$, we use the beginning of the $5$-term exact 
sequence \eqref{5 term exact} for $M={}_R k$ and $N = S$: 
\begin{equation}
\label{ex seq2 for ext k to R}
0 \lra \uext_S^1(S \otimes_R k, S) \lra 
\uext^1_R(k,S) \lra \uhom_S(\utor_1^R(S,k), S) \lra \dots 
\end{equation}

Now by Lemma~\ref{tor from S}(\ref{tor from S 1}), $\utor^R_i(S,k)$ is in $S\cat$ for 
all $i \geq 0$; in particular, $\GK_S (S \otimes_R k) \leq 1$ and 
$\GK_S(\utor^R_1(S,k)) \leq 1$.  Then  $\uext_S^1(S \otimes_R k, S) = 0$ by the 
Cohen-Macaulay property of $S$ (Lemma~\ref{CM for S}(2)) and $\uhom_S(\utor_1^R(S,k), 
S) = 0$ since $S$ is a domain with $\GK(S) = t+1 > 1$.  Thus by \eqref{ex seq2 for 
ext k to R} $\uext^1_R(k,S) = 0$, and by \eqref{ex seq for ext k to R} $\uext_R^1(k, 
R) = 0$ as well.
\end{proof}

%__________________________________________________________________________
\section{The maximal order property}
\label{ideals} \label{max order section} 

Let $A$ be a noetherian domain with Goldie quotient ring $Q$.  We say $A$ is a 
\emph{maximal order} \label{maxorder-index} in $Q$ if given any ring $T$ with $A 
\subseteq T \subseteq Q$ and nonzero elements $a, b$ of $A$ with $aTb \subseteq A$, 
we have $T = A$.  If $A$ is commutative, then $A$ is a maximal order if and only if 
$A$ is integrally closed in its fraction field \cite[Proposition 5.1.3]{MR}.  

We are interested in an equivalent formulation of the maximal order property.  For 
any left ideal $I$ of $A$, we define $\mc{O}_r(I) = \{q \in Q | Iq \subseteq I \}$ 
and $\mc{O}_l(I) = \{q \in Q | qI \subseteq I \}$.  Then $A$ is a maximal order if 
and only if $\mc{O}_r(I) = A = \mc{O}_l(I)$ for all nonzero ideals $I$ of $A$ 
\cite[Proposition~5.1.4]{MR}.  If $A$ is an $\mb{N}$-graded algebra with a graded 
ring of fractions $D$, then for any homogeneous ideal $I$ of $A$ we may also define 
$\mc{O}_r^g(I) = \{q \in D | Iq \subseteq I\}$ and $\mc{O}^g_l(I) = \{q \in D | qI 
\subseteq I\}$.  In the graded case we have the following criterion for the maximal 
order property. 
\begin{lemma}
\label{max to gr-max} Let $A$ be an $\mb{N}$-graded noetherian domain which has a 
graded quotient ring $D$ and Goldie quotient ring $Q$.  Then 
$A$ is a maximal order if and only if $\mc{O}_r^g(I) = A = \mc{O}_l^g(I)$ holds for 
all homogeneous nonzero ideals $I$ of $A$.
\end{lemma}
This result is stated in \cite[Lemma~2]{MVDB/Oyst}, but the reference given there is 
faulty and so we will supply a brief proof here.
\begin{proof}  
We may write $D \cong T[z, z^{-1}; \sigma]$ for some division ring $T$ and 
automorphism $\sigma$ of $T$. 
Then since $T$ is a maximal 
order, it follows by \cite[Propositions~IV.2.1, V.2.3]{Maury} that $D$ is a maximal 
order in $Q$. 

Assume that $\mc{O}_r^g(J) = A = \mc{O}_l^g(J)$ for all 
homogeneous ideals $J$ of $A$.  Let $I$ be any ideal of $A$, and let $q \in 
\mc{O}_r(I)$.  Then $DI$ is a 2-sided ideal of $D$ \cite[Theorem~9.20]{GW}, and also 
$q \in \mc{O}_r(DI)$.  Since $D$ is a maximal order in $Q$, this forces $q \in D$.  

Given any $d = \sum d_i \in D$ where $d_i \in D_i$, let $n$ be maximal such that $d_n 
\not = 0$ and set $\wt{d} = d_n$.  Let $\wt I$ be the 2-sided homogeneous ideal 
generated by $\wt{a}$ for  all $a \in I$.  Write $q = \sum_{i = m}^n d_i$; then since 
$Iq \subseteq I$, we have $\wt{I}d_n \subseteq \wt{I}$ and so $d_n \in \mc{O}_r^g(\wt 
I) = A$.  Then $q -d_n \in \mc{O}_r(I)$.  By induction on $n-m$ we get that $q -d_n 
\in A$ and so $q \in A$.  Thus $\mc{O}_r(I) = A$, and an analogous argument gives 
$\mc{O}_l(I) = A$, so $A$ is a maximal order.  The opposite implication is trivial.
\end{proof}

Let $S = S(\varphi)$ and $R = R(\varphi,c)$ and assume the critical density 
condition.  Our next goal is to show that $R = R(\varphi,c)$ is a maximal order.   
First, we note that the ring $S$ has this property. 
\begin{lemma}
\label{S max order}
$S = S(\varphi)$ is a maximal order. 
\end{lemma} 
\begin{proof}
By \cite[Theorem~5.11]{Zhang}, $S$ is ungraded Cohen-Macaulay and Auslander-regular, 
since $U$ has both properties; also, since $S$ is graded it is trivially stably 
free.  By \cite[Theorem~2.10]{St}, any ring satisfying these three properties is a 
maximal order.
\end{proof}

We will also require the following lemma concerning the annihilators of modules in $R\Cat$.
\begin{lemma}
\label{annihilators of modules in cat} 
\begin{enumerate}
\item
If $M \in R\Cat$, then either $_R M$ is torsion or else $\ann_R M = 0$.  
\item 
In particular, if $I$ is a nonzero ideal of $R$ then $_R(SIS/IS)$ is torsion.
\end{enumerate}
\end{lemma}
\begin{proof} 
(1) Consider the $S$-point module $\pointmod{c_i}$ for some $i \in \mb{Z}$.  By 
Lemma~\ref{truncating point modules}(1), $\pointmod{c_i}$ has point sequence $(c_i, 
c_{i-1}, c_{i-2}, \dots)$.  Then $\ann_R \pointmod{c_i} = \cap_{j = 0}^{\infty} 
\m_{c_{i - j}} \cap R$, and by the critical density of the points $\{c_i\}$ we conclude that 
$\ann_R \pointmod{c_i} = 0$.  Now the statement follows easily from the definition of $R\Cat$.

(2) Since $M ={} _R(SIS/IS)$ is a homomorphic image of a direct sum of copies of $(S/R)$, we have
$M \in R\Cat$.  Since also $IM = 0$, by part (1) $_R M$ is torsion.
\end{proof}

Recall that $R$ and $S$ have the same graded quotient ring $D$ (Lemma~\ref{quotient 
rings}).  For any graded left $R$-submodules $M,N$ of $D$, we identify $\uhom_R(M,N)$ 
with $\{d \in D \mid Md \subseteq N\}$.  Similarly, if $M,N$ are graded left 
$S$-submodules of $D$ we identify $\uhom_S(M,N)$ and $\{d \in D \mid Md \subseteq 
N\}$.  
\begin{proposition}
\label{inside S}
Let $I$ be a nonzero homogeneous ideal of $R$.  Then $\mc{O}_l^g(I) 
\subseteq S$ and $\mc{O}_r^g(I) \subseteq S$.
\end{proposition}
\begin{proof} Consider $\mc{O}_r^g(I)$ for some nonzero homogeneous ideal $I$ of $R$.  
We have that
\[
\mc{O}_r^g(I) = \{q \in D | Iq \subseteq I \} \subseteq \{q \in D |SIq \subseteq SI 
\} = \uhom_S (SI, SI).
\]
We will show that $\uhom_S (SI, SI) \subseteq S$.   Since $S$ is a maximal order by 
Lemma~\ref{S max order}, we know that $\mc{O}^g_r(SIS) = \uhom_S (SIS,SIS) = S$.  Set 
$M = SIS/SI \in S\gr$.  

Now from the exact sequence of $R$-bimodules $0 \ra SI \ra SIS \ra M \ra 0$, we get 
the following long exact sequence in $\uext$: 
\begin{gather*}
0 \ra \uhom_S(M, SIS) \ra \uhom_S(SIS, SIS) \ra \uhom_S(SI, SIS) \\
\ra \uext^1_S(M, SIS) \ra \dots
\end{gather*}
in the which the terms are again $R$-bimodules and the maps are all $R$-bimodule 
maps.  By a right sided version of Lemma~\ref{annihilators of modules in cat}(2),
which is valid by Lemma~\ref{switch sides}(1), $M_R$ must be torsion.
Then since the left $R$-structure of 
$\uext^1_S(M, SIS)$ comes from the right side of $M$, it follows from the fact that 
$_S M$ is finitely generated that the left $R$-module structure of $\uext^1_S(M, 
SIS)$ is also torsion.   

Now $\uhom_S(M,SIS) = 0$, since $S$ is a domain and $\GK(M) < \GK(S)$.  Thus 
$\uhom_S(SI, SIS)$ is an $R$-subbimodule of $D$ which is an essential left 
$R$-module extension of $\uhom_S(SIS, SIS) = S$ by a torsion module.  But by  
Lemma~\ref{ext k to R is 0}, $\uext_R^1(k, S) = 0$ and so $S$ has no nontrivial torsion extensions.  
This forces $\uhom_S(SI, SIS) = S$, and thus $\mc{O}^g_r(I) = \uhom_S(SI, SI) \subseteq \uhom_S(SI, SIS) = S$.

The proof that $\mc{O}^g_l(I) \subseteq S$ follows by applying the same argument in 
the extension of rings $R^{op} \subseteq S^{op}$, and again invoking 
Lemma~\ref{switch sides}(1).
\end{proof}

Now we may complete the proof that $R$ is a maximal order.  
\begin{theorem}
\label{max order} Assume the critical density condition, so that $R$ is noetherian.  
Then $R$ is a maximal order.
\end{theorem}

\begin{proof}
Let $I$ be any nonzero homogeneous ideal of $R$.  Then $\uhom_R(_R I,{}_R I) = 
\mc{O}_r^g(I) \subseteq S$, by Proposition~\ref{inside S}.  Set $M = 
(\uhom_R(I,I))/R$; then $_R M$ is a submodule of $_R (S/R)$, so $M \in R\Cat$.  Since 
$IM = 0$, Proposition~\ref{annihilators of modules in cat}(1) implies that $M$ is a 
torsion module.  But $\uext_R^1(k, R) = 0$ by Lemma~\ref{ext k to R is 0}, and so $R$ 
may not have any nontrivial torsion extensions.  Since $R \subseteq \uhom_R(I,I)$ is 
an essential extension, this forces $\mc{O}_r^g(I) = \uhom(I,I) = R$.  Applying the 
same argument in $R^{op}$, we get $\mc{O}_l^g(I) =R$ as well.  Thus $R$ is a maximal 
order by Lemma~\ref{max to gr-max}.
\end{proof}

%_________________________________________________________________________________
\section{The $\chi$ condition and $R\proj$}
\label{chi}  
We begin this section by reviewing some definitions from the theory of noncommutative 
projective schemes which we will use in the next two sections.  See \cite{AZ94} for 
more details.  

Let $A$ be a noetherian $\mb{N}$-graded ring which is finitely generated in degree 
one.  Let $A\Tors$ be the full subcategory of torsion objects in $A\Gr$.  Then 
$A\Tors$ is a \emph{localizing subcategory} of $A\Gr$, which means that the quotient 
category $A\Qgr = A\Gr/A\Tors$ is defined, and the exact \emph{quotient functor} 
$\pi: A\Gr \to A\Qgr$ has a right adjoint $\omega$, which is called the \emph{section 
functor}.  For torsionfree $M \in A\Gr$ we may describe $\omega \pi(M)$ explicitly as 
the unique largest essential extension $M'$ of $M$ such that $M'/M$ is torsion.  For 
all $\mc{M} \in A\qgr$, $\omega(\mc{M})$ is torsionfree and $\pi\omega(\mc{M}) \cong 
\mc{M}$. 

The \emph{noncommutative projective scheme} $A\Proj$ is defined to be the ordered 
pair $(A\Qgr, \pi(A))$, where $\pi(A)$ is called the \emph{distinguished object}.  We 
write $A\Proj \cong B\Proj$ if there is an equivalence of categories $A\Qgr \sim 
B\Qgr$ under which the distinguished objects correspond.  We also work with the 
subcategories of noetherian objects $A\gr, A\tors, A\qgr = A\gr/A\tors$, and we set $A\proj = 
(A\qgr, \pi(A))$.

The category $A\Qgr$ has enough injectives and so $\ext$ groups are defined in this 
category.  The shift functor $M \mapsto M[1]$, which is an autoequivalence of the 
category $A\Gr$, descends naturally to an autoequivalence of $A\Qgr$.  For $\mc{M}, 
\mc{N} \in A\Qgr$ we define 
\[
\uext^i(\mc{M},\mc{N}) = \oplus_{i = - \infty}^{\infty} \ext^i_{A\Qgr}(\mc{M}, 
\mc{N}[i]).
\]
Now we define cohomology and graded cohomology for $A\Proj$ by setting $\coH^i(\mc{N}) 
= \ext^i(\pi(A), \mc{N})$ and $\ucoH^i(\mc{N}) = \uext^i(\pi(A), \mc{N})$ for  
$\mc{N} \in A\Qgr$.  The section functor $\omega$ may be described using cohomology 
as $\omega(\mc{M}) = \ucoH^0(\mc{M})$ for all $\mc{M} \in A\Qgr$.  

For making explicit computations, it is useful to note that for all $M, N \in A\gr$ and all $p \geq 0$,
\begin{equation}
\label{cohom form 1}
\uext^p_{A\Qgr}(\pi(M), \pi(N)) \cong \lm \uext^p_A(M_{\geq n}, N).
\end{equation}
In case  $M = A$, we have the following additional formula for all $p \geq 1$ (see \cite[Proposition~7.2(2)]{AZ94}):
\begin{equation}
\label{cohom form 2}
\ucoH^p(\pi(N)) \cong \lm \uext^{p+1}_A(A/A_{\geq n}, N).
\end{equation}
If $A$ is commutative this last formula amounts to the usual 
correspondence between sheaf cohomology in $A\proj$ and local cohomology for the ring 
$A$.

Recall the $\chi$ conditions which were defined in \S\ref{homological lemmas}. Artin 
and Zhang proved the following noncommutative analog of Serre's theorem:
\begin{theorem} \cite[Theorem~4.5]{AZ94}
\label{noncomm Serre}
Let $A$ be a right noetherian $\mb{N}$-graded algebra, and let $\mc{A} = \pi(A)$ be 
the distinguished object of $A\proj$.  Then $B = \bigoplus_{i \geq 0} 
\coH^0(\mc{A}[i])$ is naturally a graded ring and there is a canonical homomorphism 
$\psi: A \ra B$.  If $A$ satisfies $\chi_1$ then $\psi$ is an isomorphism in large 
degree, and $B\proj \cong A\proj$. \hfill $\Box$
\end{theorem}
In other words, if $A$ satisfies $\chi_1$ then the noncommutative projective scheme 
$A\proj$, together with the shift functor $M \mapsto M[1]$, determines the ring $A$ 
up to a finite dimensional vector space.

In this section, we will analyze the $\chi$ conditions for $R = R(\varphi,c)$, 
assuming the critical density condition throughout.  We shall show that $R$ satisfies 
$\chi_1$, but that $\chi_i$ fails for $R$ for $i \geq 2$.  In particular, $R$ 
satisfies the noncommutative Serre's theorem, and is the first example of a 
noetherian algebra which satisfies $\chi_1$ but not all of the $\chi$ conditions.

We will say that $\chi_i(M)$ \emph{holds} for a particular module $M \in A\Gr$ if 
$\uext^j_A(k, M) < \infty$ for $0 \leq j \leq i$.  The reader may easily prove the 
following simple facts.
\begin{lemma}
\label{reducing exact sequences}
Let $0 \ra M' \ra M \ra M'' \ra 0$ be an exact sequence in $R\gr$, and let $N \in R\gr$.  
\begin{enumerate} 
\item \label{reducing exact sequences 1} 
If $\chi_1(M')$ and $\chi_1(M'')$ hold then $\chi_1(M)$ holds. 
\item \label{reducing exact sequences 2} 
If $\chi_1(M)$ holds then $\chi_1(M')$ holds. 
\item \label{reducing exact sequences 3} 
If $\dim_k N < \infty$ then $\chi_1(N)$ holds. \hfill $\Box$ 
\end{enumerate}
\end{lemma}

To prove $\chi_1$ for $R$ we will reduce to the case of $S$-modules. 
\begin{proposition}
\label{chi for S modules} Suppose that $N \in S\gr$.  Then 
$\chi_1(_R N)$ holds.
\end{proposition}
\begin{proof}
Consider the first 3 terms of the 5-term exact sequence \eqref{5 term exact} for $M = 
{}_R k$:
\[
0 \lra \uext_S^1(S \otimes_R k, N) \lra \uext^1_R(k,N) \lra \uhom_S(\utor_1^R(S,k),
N) \lra \dots
\]
Now $\utor_i^R(S,k)$ is in $S\cat$ for any $i \geq 0$, by Lemma~\ref{tor from 
S}(\ref{tor from S 1}).  Then by Proposition~\ref{exts over S}, we conclude that 
$\dim_k \uext_S^j(\utor_i^R(S,k)), N) < \infty$ for $j = 0,1$ and $i \geq 0$.  Thus 
$\dim_k \uext^1_R(k,N) < \infty$.
\end{proof}

The following result completes the proofs of Theorems~\ref{main result 2} and 
\ref{main result 3} from the introduction.
\begin{theorem}
\label{chi for R} \label{chi fails} Assume the critical density condition and let $R = 
R(\varphi,c)$. 
\begin{enumerate}
\item \label{chi for R 1} 
$R$ satisfies $\chi_1$. 
\item \label{chi for R 2}  
$\uext_R^2(k, R)$ is not bounded, and $\chi_i$ fails for all 
$i \geq 2$. 
\end{enumerate}
\end{theorem}
\begin{proof}  
(1) By Lemma~\ref{reducing exact sequences}(1) and induction it is enough to show 
that $\chi_1(M)$ holds for all graded cyclic $R$-modules $M$. 

Let $R/I$ be an arbitrary graded cyclic left $R$-module.  If $I = 0$, then $\chi_1(_R 
R)$ holds by Lemma~\ref{ext k to R is 0}.  Assume then that $I \not = 0$.  Consider 
the exact sequence \eqref{exact seq}.  Now $\chi_1(_R(S/SI))$ holds by 
Proposition~\ref{chi for S modules}.  By Lemma~\ref{comparing cyclic R and S 
modules}(1), both $(SI \cap R)/I$ and $S/(R + SI)$ have finite filtrations with 
factors which are either torsion or shifted $R$-point modules with a compatible 
$S$-module structure.  Then $\chi_1((SI \cap R)/I)$ and $\chi_1(S/(R + SI))$ hold, by 
Proposition~\ref{chi for S modules} and Lemma~\ref{reducing exact 
sequences}(1),(3).   Finally, applying Lemma~\ref{reducing exact sequences}(1),(2) to 
\eqref{exact seq} we get that $\chi_1(R/I)$ holds.

(2) Consider the long exact sequence in $\uext_R(k, -)$ that arises from the short
exact sequence of $R$-modules $0 \ra R \ra S \ra S/R \ra 0$:
\begin{equation}
\label{eq chi fails}
\dots \ra  \uext_R^1(k, S) \ra \uext_R^1(k, S/R) \ra \uext_R^2(k, R) \ra \dots
\end{equation}
Now $\uext_R^1(k, S) = 0$, by Lemma~\ref{ext k to R is 0}.  On the other hand,
\[
\uext_R^1(k, S/R) \cong \oplus_{i=1}^{\infty} \uext_R^1(k, \pointmod{c_{-1}})[-i]
\]
by Corollary~\ref{structure of S/R}(2), since $\uext$ commutes with direct sums in 
the second coordinate \cite[Proposition~3.1(1)(b)]{AZ94}.  By 
Theorem~\ref{classifying point ideals}(2), it is clear that the point module 
$\pointmod{c_{-1}}$ has a nontrivial extension by $k[1]$, since any point module 
$P(c_0, e)$ defined there satisfies $(P(c_0,e)[1])_{\geq 0} \cong{} _R 
\pointmod{c_{-1}}$. Thus $\uext_R^1(k, \pointmod{c_{-1}}) \neq 0$, and so 
$\oplus_{i=1}^{\infty} \uext_R^1(k, \pointmod{c_{-1}})[-i] \cong \uext_R^1(k, S/R)$ 
is not right bounded.  Then by the exact sequence \eqref{eq chi fails}, $\uext_R^2(k, 
R)$ is also not right bounded.  In particular, $\dim_k \uext_R^2(k,R) = \infty$ and 
$\chi_i$ fails for $R$ for all $i \geq 2$ by definition.
\end{proof}

We see next that the failure of $\chi_i$ for $R$ for $i \geq 2$ is reflected in the 
cohomology of $R\proj$.  We recall the noncommutative version of Serre's finiteness 
theorem which was proved by Artin and Zhang, which we have restated slightly.
\begin{theorem} \cite[Theorem~7.4]{AZ94}
\label{Serres finiteness} Let $A$ be a left noetherian finitely $\mb{N}$-graded 
algebra which satisfies $\chi_1$.  Then $A$ satisifes $\chi_i$ for some $i \geq 2$ if 
and only if the following two conditions hold: 
\begin{enumerate}
\item \label{Serres finiteness 1}
$\dim_k \coH^j(\mc{N}) < \infty$ for all $0 \leq j < i$ and all $\mc{N} \in A\qgr$. 
\item \label{Serres finiteness 2}
$\ucoH^j(\mc{N})$ is right bounded for all $1 \leq j < i$ and all $\mc{N} \in A\qgr$. 
\end{enumerate}
\end{theorem}
\begin{proof}
This follows immediately from the proof of \cite[Theorem~7.4]{AZ94}.
\end{proof}

\begin{lemma}
\label{chi and ext-finiteness} Let $A$ be a left noetherian finitely $\mb{N}$-graded 
algebra satisfying $\chi_i$.  Then $\dim_k \ext^j(\mc{M}, \mc{N}) < \infty$ for $0 
\leq j < i$ and for all $\mc{M}, \mc{N} \in A\qgr$.
\end{lemma}
\begin{proof}
Let $\mc{A} = \pi(A)$.  Since any $M \in A\gr$ is an image of some finite sum of 
shifts of $A$, in $A\qgr$ there is an exact sequence
\[
0 \ra \mc{M}' \ra \mc{F} \ra \mc{M} \ra 0
\]
where we have $\mc{F} =  \oplus_{i = 1}^n \mc{A}[d_i]$ for some integers $d_i \in \mb{Z}$.  
Then $\ext^j(\mc{F}, \mc{N}) \cong \oplus_{i = 1}^n \ext^j(\mc{A}, \mc{N}[-d_i]) = \oplus_{i =1}^n 
\coH^j(\mc{N}[-d_i])$ and so $\dim_k \ext^j(\mc{F}, \mc{N}) < \infty$ for all $0 \leq 
j < i$ by Theorem~\ref{Serres finiteness}.

We induct on $j$.  If $j = 0$ then there is an exact sequence $0 \ra \hom(\mc{M}, 
\mc{N}) \ra \hom(\mc{F}, \mc{N})$ from which it follows that $\dim_k \hom(\mc{M}, 
\mc{N}) < \infty$. For $0 < j < i$, there is the long exact sequence
\[
\dots \ra \ext^{j-1}(\mc{M}', \mc{N}) \ra \ext^j(\mc{M}, \mc{N}) \ra \ext^j(\mc{F}, 
\mc{N}) \ra \dots 
\] 
and since $\dim_k \ext^{j-1}(\mc{M}', \mc{N}) < \infty$ by the induction hypothesis, 
we have $\dim_k \ext^j(\mc{M}, \mc{N}) < \infty$ as well.  This completes the induction step 
and the proof.
\end{proof}

We can now make the failure of the Serre's finiteness theorem for $R\proj$ explicit.
\begin{lemma}
\label{R not ext-finite} 
Let $\mc{R} = \pi(R) \in R\qgr$ be the distinguished object 
of $R\proj$.  Then $\dim_k \coH^1(\mc{R}) = \infty$.
\end{lemma}
\begin{proof}
Set $\mc{S} = \pi(S) \in R\Qgr$.  The exact sequence $0 \ra R \ra S \ra S/R \ra 0$ 
descends to an exact sequence $0 \ra \mc{R} \ra \mc{S} \ra \mc{S}/\mc{R} \ra 0$ in 
$R\qgr$.  For $\mc{M} \in R\Qgr$, the cohomology $\coH^0(\mc{M})$ may be identified 
with the zeroeth graded piece of the module $\omega(\mc{M})$, where $\omega$ is the 
section functor.  Recall also that for torsionfree $M \in A\Gr$, $\omega\pi(M)$ is 
the largest essential extension of $M$ by a torsion module.  Since $\uext^1_R(k, S) = 
0$ by Lemma~\ref{ext k to R is 0}, $_R S$ has no nontrivial torsion extensions and so 
$\omega(\mc{S}) = S$.  In particular, $\dim_k \coH^0(\mc{S}) = \dim_k S_0 = 1$. On 
the other hand,  $S/R = \oplus_{i =1}^{\infty} \pointmod{c_{-1}}[-i]$ is an infinite 
direct sum of shifted $R$-point modules by Corollary~\ref{structure of S/R}(2).  For 
each $i \geq 0$, by Lemma~\ref{classifying point ideals}(2) there is some $R$-point 
module $P(c_i, e_i)$ which satisfies $P(c_i, e_i)_{\geq i+1} \cong 
\pointmod{c_{-1}}[-i-1]$.  Then $M = \oplus_{i =0}^{\infty} P(c_i, e_i)$ is an 
essential extension of $S/R$ by a torsion module, so $M \subseteq 
\omega(\mc{S}/\mc{R})$ and it follows that $\dim_k \coH^0(\mc{S}/\mc{R}) \geq \dim_k 
M_0 = \infty$.  Now the long exact sequence in cohomology forces $\dim_k 
\coH^1(\mc{R}) = \infty$ as well.
\end{proof}

The following result, which proves Theorem~\ref{main result 4} from the introduction, 
shows that the category $R\qgr$ is necessarily quite different from any of the 
standard examples.  
\begin{theorem}
\label{R proj not comm} Assume the critical density condition. 
\begin{enumerate} 
\item \label{R proj not comm 1} 
Suppose that $A$ is a left noetherian finitely $\mb{N}$-graded $k$-algebra which 
satisfies $\chi_2$. Then the categories $A\qgr$ and $R\qgr$ are not equivalent. 
\item \label{R proj not comm 2} 
$R\qgr$ is not equivalent to $\coh X$, the category of coherent sheaves on $X$, 
for any commutative projective scheme $X$.
\end{enumerate}
\end{theorem}
\begin{proof}
(1) The proof is immediate from Lemmas~\ref{chi and ext-finiteness} and \ref{R not 
ext-finite}.  

(2) This follows from part (1) and the usual commutative Serre's theorem.  
\end{proof}
%______________________________________________________________________________
\section{Global and cohomological dimension of $R\proj$}
\label{cohom dim} In this section, our goal 
is to show that $R\proj$ has finite global dimension, and thus finite cohomological 
dimension.  We will also give upper bounds for these numbers. 

Let us recall the definitions of these concepts: 
\begin{definition}
\label{cohomdim-index} Let $A$ be a connected finitely generated $\mb{N}$-graded 
algebra.  The \emph{global dimension} of $A\qgr$ (or $A\proj$) is 
\[
\gld(A\qgr) = \sup \{i \mid \ext^i(\mc{M}, \mc{N}) \neq 0\ \text{for some}\ \mc{M}, 
\mc{N} \in A\qgr \}. 
\]
The \emph{cohomological dimension} of $A\proj$ is  
\[
\cd(A\proj) = \sup \{i \mid \coH^i(\mc{N}) \neq 0\ \text{for some}\ \mc{N} \in A\qgr 
\}. 
\]
\end{definition}
If $A\qgr$ has finite global dimension, then it is immediate that $A\proj$ has finite 
cohomological dimension.  We remark that it is not known if there exists any graded 
algebra $A$ such that $\cd(A\proj) = \infty$.

Now let $S = S(\varphi)$ and $R = R(\varphi,c)$, and assume the critical density 
condition as usual.  It is easy to compute the global and cohomological dimensions of 
$S\proj$:
\begin{lemma}
\label{cd S} $\cd(S\proj) = \gld(S\qgr) = \GK(S) -1 = t$.  
\end{lemma}
\begin{proof}
Since $S$ is a Zhang twist of the polynomial ring $U$, we have an isomorphism $S\proj 
\cong (\coh \mb{P}^t, \mc{O}_{\mb{P}^t})$, and the values of both dimensions for the 
commutative scheme $\mb{P}^t$ are well known. 
\end{proof}

The main machinery we will use to study $\ext$ groups in $R\qgr$ is 
Proposition~\ref{spectral sequence collapses} below, 
which needs the spectral sequence given in the following lemma.  
\begin{lemma}
\label{lim spec seq} For any system in $R\Gr$ of the form $\dots \ra M_n \ra \dots 
\ra M_1 \ra M_0$ and any $N \in S\Gr$ there is a convergent spectral sequence  of the 
form
\[
E^{pq}_2 = \lm \uext^p_S(\utor^R_q(S, M_n), N) \underset{p}{\Rightarrow} \lm 
\uext^{p+q}_R(M_n,N). 
\]
\end{lemma}
\begin{proof}
Consider the spectral sequence \eqref{spectral sequence} for arbitrary $_R M \in 
R\Gr$:
\[ 
\uext^p_S(\utor^R_q(S, M), N) \underset{p}{\Rightarrow} \uext^{p+q}_R(M,N).
\]
Let $\mc{C}$ be the category of all $\mb{N}$-indexed directed systems of modules in 
$R\Gr$ of the form 
\[
\dots \ra M_n \ra \dots \ra M_1 \ra M_0. 
\]
Let $\mc{D}$ be the analogous category of directed systems of modules in $S\Gr$.  
Both of these categories have enough projectives and injectives.  For example, if $P$ 
is a projective object of $R\Gr$, then any object in $\mc{C}$ of the form
\begin{equation}
\label{proj objects} \dots \ra 0 \ra 0 \ra 0 \ra P \overset{\cong}{\ra} P 
\overset{\cong}{\ra} \dots \overset{\cong}{\ra} P
\end{equation}
is projective, and clearly every object in $\mc{C}$ is an image of a direct sum of 
objects of this form.  See \cite[Exercises 2.3.7, 2.3.8]{We} for more details. The 
functor $S \otimes_R - :R\Gr \ra S\Gr$ extends to a functor $G: \mc{C} \ra \mc{D}$.  
We also have a functor $F: \mc{D} \ra \operatorname{Ab}$ defined by $\{ L_n \}_{n \in 
\mb{N}} \mapsto \lm \uhom_S(L_n, N)$, where $\operatorname{Ab}$ is the category of 
abelian groups.  It is easy to see that $G$ is right exact; since $\operatorname{Ab}$
has exact direct limits \cite[Theorem 2.6.15]{We}, $F$ is left exact. Finally, $G$ 
sends any direct sum of objects in $\mc{C}$ of the form in \eqref{proj objects} to a 
projective object in $\mc{D}$.  Then corresponding to the composition of functors $F 
\circ G$ is a Grothendieck spectral sequence (see \cite[Theorem 11.40]{Rotman})
\[
E^{p,q}_2 = R^pF(L_qG(M_{.})) \underset{p}{\Rightarrow} R^{p+q}(FG)(M_{.})
\]
which we leave to the reader to show unravels to the spectral sequence 
required by the lemma.
\end{proof}

To get the most out of the spectral sequence, we also note the following simple lemma.
\begin{lemma}
\label{ignore torsion} 
Let $A$ be a connected graded noetherian ring, and let $\dots 
\ra M_n \ra \dots \ra M_1 \ra M_0$ a directed system of modules in $A\gr$.  For each 
$n$, let $\tau(M_n)$ be the torsion submodule of $M_n$.  Then for any $N \in A\gr$ 
and $p \geq 0$, 
\[
\lm \uext^p_A(M_n, N) \cong \lm \uext^p_A(M_n/\tau(M_n), N)
\]
as $k$-spaces.
\end{lemma}
\begin{proof} 
Since direct limits are exact in the category of abelian groups \cite[Theorem 
2.6.15]{We}, there is a long exact sequence in $\lm \uext(-, N)$ arising from short 
exact sequence of complexes $0 \to \tau(M.) \to M. \to M./\tau(M.) \to 0$.    Given 
fixed $n$, the module $\tau(M_n)$ is bounded and so for some $n' \gg n$ the natural 
map $\tau(M_{n'}) \ra \tau(M_n)$ is zero, and then the natural map 
$\uext^p_A(\tau(M_{n}), N) \ra \uext^p_A(\tau(M_{n'}), N)$ is zero.  Thus $\lm 
\uext^p_A(\tau(M_n), N) = 0$ for all $p \geq 0$, and the desired result follows from 
the long exact sequence.
\end{proof}

\begin{proposition}
\label{spectral sequence collapses} \label{relation between R and S cohomology}
 Let $N \in S\Gr$, and let $M \in R\gr$. 
\begin{enumerate}
\item \label{spectral sequence collapses 1} 
As graded vector spaces, for all $m \geq 0$ we have 
\[
\uext_{R\Qgr}^{m}(\pi(M), \pi(N)) \cong \lm \uext^m_S(S \otimes_R M_{\geq n}, N).
\]
\item \label{spectral sequence collapses 2} 
In case $M = R$, for $m \geq 1$ we have 
\[
\ucoH^m(\pi(N)) \cong \lm \uext^{m+1}_S(S/SR_{\geq n}, N).
\]
\end{enumerate}
\end{proposition}
\begin{proof}
(1) We use the spectral sequence of Lemma~\ref{lim spec seq}:
\[
E^{pq}_2 = \lm \uext^p_S(\utor^R_q(S, M_{\geq n}), N) \underset{p}{\Rightarrow} \lm 
\uext^{p+q}_R(M_{\geq n},N). 
\]
Our goal is to show that $E^{pq}_2 = 0$ for any pair of indices $p,q$ with $q \geq 
1$.  

Fix $q \geq 1$.  For fixed $n \geq 0$, we claim first that there is some $n' \geq n$ 
such that the natural map $\psi_1: \utor^R_q(S, M_{\geq n'}) \ra \utor^R_q(S, M_{\geq 
n})$ is $0$.  As in Lemma~\ref{tor from S}, there is a right point module $Q$ of $R$ 
such that $S/R \cong \bigoplus_{i = 1}^{\infty} Q[-i]$ as right $R$-modules.  Now by 
the fact that $\utor$ commutes with direct sums \cite[Proposition~2.4(1)]{AZ94} we 
get a commutative diagram  
\[
\begin{array}{ccc}
\utor_q^R (S, M_{\geq n'}) & \overset{\psi_1}{\lra}  & \utor_q^R (S,M_{\geq n}) \\
\big\da & &  \big\da \\
\utor_q^R (S/R,M_{\geq n'} ) & \overset{\psi_2}{\lra} & \utor_q^R (S/R,M_{\geq n}) \\
\big\da\cong & &  \big\da\cong \\
\bigoplus_{i=1}^{\infty} \utor_q^R (Q,M_{\geq n'})[-i] & \overset{\psi_3}{\lra} & 
\bigoplus_{i=1}^{\infty} \utor_q^R (Q,M_{\geq n})[-i]
\end{array}
\]
where the top two vertical maps are at least injections (see Lemma~\ref{tor from S})
and the $\psi_i$ are the natural maps.  

Now $T^n = \utor_q^R(Q, M_{\geq n})$ is bounded, by Lemma~\ref{tor from S}(2).  Also, 
clearly the left bound $l(n)$ of $T^n$ satisfies $\lm l(n) = \infty$.  It follows 
that for $n' \gg n$ the natural map $\theta: T^{n'} \ra T^n$ is $0$.  The restriction 
of the map $\psi_3$ to any summand is just a shift of the map $\theta$, so $\psi_3 = 
0$ for $n' \gg n$.  Finally, the commutative diagram gives $\psi_1 = 0$ for $n' \gg 
n$.  This proves the claim.

Write $_n E_2^{pq} = \uext^p_S(\utor^R_q(S, M_{\geq n}), N)$.  Since $\psi_1 = 0$ for 
$n' \gg n$, the natural map ${}_n E^{pq}_2 \ra {}_{n'} E^{pq}_2$ is also zero for $n' 
\gg n$.  Since $n$ was arbitrary, we have $E^{pq}_2 = \lm {}_n E^{pq}_2 = 0$.

Therefore only the $E^{pq}_2$ with $q=0$ are possibly nonzero, and the spectral 
sequence collapses, giving an isomorphism of vector spaces for all $m \geq 1$ as 
follows (using also \eqref{cohom form 1}):  
\begin{gather*}
\uext_{R\Qgr}^m(\pi(M),\pi(N)) \cong \lm \uext^m_R(M_{\geq n},N) \cong \lm 
\uext^m_S(\utor^R_0(S, M_{\geq n}), N) 
\\ = \lm \uext^m_S(S \otimes_R M_{\geq n}, N).
\end{gather*}
(2) Assume that $M = R$.  Setting $T^n = \utor_q^R(Q, R/R_{\geq n})$ in this case (where $Q$ is as in part (1)) 
it is obvious that $T^n$ is bounded, since $R/R_{\geq n}$ is, and it is still true for 
$q \geq 1$ that the left bound $l(n)$ of $T^n$ satisfies $\lm l(n) = \infty$ 
\cite[Proposition 2.4(6)]{AZ94}.  Then the same argument as in part (1) shows that 
the spectral sequence 
\[
E^{pq}_2 = \lm \uext^p_S(\utor^R_q(S, R/R_{\geq n}), N) \underset{p}{\Rightarrow} \lm 
\uext^{p+q}_R(R/R_{\geq n},N)
\]
also collapses, so we have (using \eqref{cohom form 2}) that
\begin{gather*}
\ucoH^m_R(\pi(N)) \cong \lm \uext^{m+1}_R(R/R_{\geq n},N) \cong \lm 
\uext^{m+1}_S(\utor^R_0(S, R/R_{\geq n}), N) 
\\ = \lm \uext^{m+1}_S(S/SR_{\geq n}, N)
\end{gather*}
for all $m \geq 1$.
\end{proof}
Armed with the preceding results, we can now show that $R\qgr$ has 
finite global dimension, and calculate upper bounds on the values of the global 
dimension and cohomological dimension for $R\proj$.  This proves Theorem~\ref{main 
result 5} from the introduction.
\begin{theorem}  
\label{cohom dim R} Assume the critical density condition for $R = R(\varphi,c)$, and 
recall that $t = \GK(R) -1$.  Then 
\begin{enumerate} 
\item \label{cohom dim R 1}
$\gld(R\qgr) \leq t+1$. 
\item \label{cohom dim R 2} 
$\cd(R\proj) \leq t$. 
\end{enumerate} 
\end{theorem}
\begin{proof}
(1) Given any $M \in R\gr$ and $N \in S\Gr$, Proposition~\ref{spectral sequence 
collapses}(\ref{spectral sequence collapses 1}) and Lemma~\ref{ignore torsion} show 
that 
\[
\uext^p_{R\qgr}(\pi(M), \pi(N)) \cong \lm \uext^p_S(S \otimes M_{\geq n}, N) \cong   
\lm \uext^p_S(N'_n, N)
\]
where $N'_n = (S \otimes M_{\geq n})/\tau(S \otimes M_{\geq n})$.  For all $n \geq 
0$, the $S$-module $N'_n$ is torsionfree.  By the graded Auslander-Buchsbaum formula 
\cite[Exercise 19.8]{E}, a torsionfree module in $U\gr$ has projective dimension at 
most $(\dim U) -1 = t$, since it has depth $\geq 1$; by the equivalence of categories 
$U\Gr \sim S\Gr$ we infer that each $N'_n$ has projective dimension at most $t$ over $S$.  
Then for $p > t$ each term of the direct limit is $0$ and so $\uext^p_{R\qgr}(\pi(M), 
\pi(N)) = 0$.

Now if $L$ is an arbitrary graded cyclic $R$-module, then either $L = R$, in which 
case we have the exact sequence $0 \ra R \ra S \ra S/R \ra 0$, or $L = R/I$ for $I 
\neq 0$, in which there is the exact sequence \eqref{exact seq}: $0 \ra (SI \cap R)/I 
\ra R/I \ra S/SI \ra S/(R + SI) \ra 0$.  Since $S/R, (SI \cap R)/I,$ and $S/(R + SI)$ 
all have finite filtrations with factors which have an $S$-structure compatible with 
their left $R$-structure (Lemmas~\ref{structure of S/R}, \ref{comparing cyclic R and S modules}), we conclude 
that $\uext^p_{R\qgr}(\pi(M), \pi(L)) = 0$ for $p > (t+1)$.  Since any $L' \in R\gr$ 
has a finite filtration by cyclic modules, we see that $\uext^p_{R\qgr}(\pi(M), 
\pi(L')) = 0$ for all $M, L' \in R\gr$ and $p > (t+1)$, and so $\gld(R\qgr) \leq 
t+1$.    

(2) We will show that for any $N \in S\Gr$, $\uext^p_{R\qgr}(\pi(R), \pi(N)) = 0$ for 
$p > (t-1)$.  Then the bound $\cd(R\proj) \leq t$ will follow by a similar argument 
as in part (1).  Let $J^{(n)}$ be the left $S$-ideal $\m_{c_0} \cap \m_{c_1} \cap 
\dots \cap \m_{c_{n-1}}$.  By Theorem~\ref{char of R}, $R_n = (J^{(n)})_n$, and 
furthermore $J^{(n)}$ is generated in degrees $\leq n$, by Lemma~\ref{products versus 
intersections}.  It is easy to see then that $J^{(n)}$ is the largest extension of 
$SR_{\geq n}$ inside $S$ by a torsion module.
 
Now by Proposition~\ref{spectral sequence collapses}(\ref{spectral sequence collapses 
2}) and Lemma~\ref{ignore torsion}, we have 
\[
\uext^p_{R\qgr}(\pi(R), \pi(N)) \cong \lm \uext^{p+1}_S(S/SR_{\geq n}, N) \cong \lm 
\uext^{p+1}_S(S/J^{(n)}, N).
\]
But each $S/J^{(n)}$ is torsionfree and hence has projective dimension at most $t$, 
by the same argument as in part (1), so every term in the direct limit is zero when 
$p 
> (t-1)$.
\end{proof}

Before leaving the subject of cohomological dimension, we wish to mention another 
approach to cohomology for noncommutative graded algebras which is provided by the 
work of Van Oystaeyen and Willaert on schematic algebras \cite{VOyWi1, VOyWi2, 
VOyWi3}.  An algebra graded $A$ is called \emph{schematic} \label{schematic-index} if 
it has enough Ore sets to give an open cover of $A\proj$; we shall not concern 
ourselves here with the formal definition.  For such algebras one can define a 
noncommutative version of \v{C}ech cohomology which gives the same cohomology groups 
as the cohomology theory we studied above.  

It turns out that the theory of schematic algebras is of no help in computing the 
cohomology of $R\proj$.  Indeed, if $A$ is a connected $\mb{N}$-graded noetherian 
schematic algebra then $\uext^n_A(_A k,A)$ is torsion as a right $A$-module for all 
$n \in \mb{N}$ \cite[Proposition~3]{VOyWi3}, hence finite dimensional over $k$.  But 
we saw in Theorem~\ref{chi fails} that $\dim_k \uext^2_R(_R k,R) = \infty$.  Thus we 
have incidentally proven the following proposition.
\begin{proposition}
\label{not schematic} Assume the critical density condition.  Then $R = R(\varphi,c)$ 
is a connected $\mb{N}$-graded noetherian domain, generated in degree one, which is 
not schematic. \hfill $\Box$
\end{proposition}
The previously known non-schematic algebras have not been generated in degree one 
\cite[page 12]{VOyWi3}. 
%_________________________________________________________________________________
\section{The critical density property}
\label{examples}

We saw in Theorem~\ref{R noetherian} that the noetherian property for $R(\varphi, c)$ 
depends on the critical density of the set of points $\mc{C} = \{ \varphi^{i}(c) 
\}_{i \in \mb{Z}}$, and we have been assuming that $\mc{C}$ is critically dense ever 
since.  In this section, we justify this assumption by showing that the critical 
density of $\mc{C}$ holds for generic choices of $\varphi$ and $c$.  We will only 
concern ourselves with pairs $(\varphi,c)$ such that Hypothesis~\ref{distinct} holds, 
that is such that $c$ has infinite order under $\varphi$.  Next, we will show that in 
case $\cha k = 0$, the set $\mc{C}$ is critically dense if and only if it is dense, 
which is a much simpler condition to check.  Finally, we discuss rings generated by 
Eulerian derivatives, which was the context in which rings of the form $R(\varphi,c)$ 
first appeared in the literature \cite{Jordan}. We translate our earlier results into 
this language, and show that they solve several open questions in \cite{Jordan}. 

Throughout this section we write $c_i = \varphi^{-i}(c)$.  As we did earlier in 
\S\ref{explicit ring}, we think of automorphisms of $\mb{P}^t$ as elements of 
$\pgl_{t+1}(k) = \gl_{t+1}(k) /k^{\times}$, which act on the left on column vectors 
of homogeneous coordinates for $\mb{P}^t$.  We write $\diag(p_0, p_1, \dots, p_t)$ 
for the automorphism which is represented by a diagonal matrix with diagonal entries 
$p_0, p_1, \dots, p_t$.

Let us define precisely our (somewhat nonstandard) intended meaning of the word 
``generic''.
\begin{definition}
\label{generic-index} \label{generic} A subset $U$ of a variety $X$ is \emph{generic} 
if its complement is contained in a countable union of proper closed subvarieties of 
$X$. 
\end{definition}
If the base field $k$ is uncountable, a generic subset is intuitively very large.  
For example, if $k = \mb{C}$ then a property which holds generically holds ``almost 
everywhere'' in the sense of Lebesgue measure.  For any results below which involve 
genericity we will assume that $k$ is uncountable.

In the next theorem we will prove that $\mc{C}$ is critically dense for $\varphi$ a 
suitably general diagonal matrix, and $c$ chosen from an open set of $\mb{P}^t$.  The 
proof will depend on the following combinatorial lemma.
\begin{lemma}
\label{determinant} Fix $d \geq 1$, and set $N = \binom{t+d}{d}$.  Let $U = k[x_0, 
x_1, \dots, x_t]$ be the polynomial ring, and give monomials in $U$ the lexicographic 
order with respect to some fixed ordering of the variables.  Let $f_1, f_2, \dots, 
f_N$ be the monomials of degree $d$ in $U = k[x_0, x_1, \dots, x_t]$, enumerated so 
that $f_1 < f_2 < \dots f_N$ in the lex order.  Fix some sequence of distinct 
nonnegative integers $a_1 < a_2 < \dots < a_N$.  Then the polynomial $\det 
(f_i^{a_j}) \in U$ is nonzero. 
\end{lemma}
\begin{proof}
Set $F = \det (f_i^{a_j}) \in U$.  Let $S_N$ be the symmetric group on $N$ elements, 
with identity element $1$; then $F$ is a sum of terms of the form $h_{\sigma} = \pm 
\prod_{i = 1}^{N} f_i^{a_{\sigma(i)}}$ for $\sigma \in S_N$.  It is straightforward 
to check that the monomial $f_1^{a_1} f_2^{a_2} \dots f_N^{a_N}$ is the unique 
largest in the lex order occurring among the $h_{\sigma}$, and that it occurs only in 
$h_1$ and thus may not be cancelled by any other term.  
\end{proof}

The next theorem includes, in particular, the result of Proposition~\ref{diag case 1} 
which we stated earlier without proof.
\begin{theorem}
\label{diag case}  Let $\varphi = \diag(1, p_1, p_2, \dots, p_t)$  with $\{p_1, p_2, 
\dots p_t \}$ algebraically independent over the prime subfield of $k$.  Let $c = 
(b_0:b_1: \dots :b_t) \in \mb{P}^t$ with $b_i \neq 0$ for all $0 \leq i \leq t$.   
Then $\mc{C} = \{ \varphi^i(c) \}_{i \in \mb{Z}}$ is critically dense and 
$R(\varphi,c)$ is noetherian.  
\end{theorem}
\begin{proof}  
We have the explicit formula $c_{-n} = (b_0 :b_1 p_1^n : b_2 p_2^n : \dots : b_t 
p_t^n)$.  We will actually prove that the set of points $\mc{C}$ is in \emph{general 
position}: that is, that at most $\binom{n+d}{d}$ of the points lie on any degree $d$ 
hypersurface of $\mb{P}^t$.  This will obviously imply that $\mc{C}$ is critically 
dense.  

Suppose that $\mc{C}$ fails to be in general position.  Then there is some $d \geq 1$ 
and a sequence of $N = \binom{t+d}{d}$ integers $a_1 < a_2 < \dots < a_N$ such that 
the points $c_{a_1}, c_{a_2}, \dots, c_{a_N}$ lie on a degree $d$ hypersurface in 
$\mb{P}^t$.  We may assume that the $a_i$ are nonnegative, since if the $\{c_{a_i}\}$ 
lie on a degree $d$ hypersurface then the same is true of the points 
$\{\varphi^{-m}(c_{a_i})\} = \{c_{a_{i}+m}\}$ for any $m \in \mb{Z}$.  Let $f_1, f_2, 
\dots f_N$ be the distinct degree $d$ monomials in the variables $x_i$ of $U$.  It 
follows that $\det (f_i(c_{a_j})) = 0$.

Given the explicit formula for $c_n$, one calculates that  
\[
\det (f_i(c_{a_j})) = B[\det (f_i^{a_j})](1:p_1^{-1}:p_2^{-1}: \dots: p_t^{-1}) = 0
\]
where $B$ is a monomial in the $b_i$ and hence is nonzero by hypothesis.  Now by 
Lemma~\ref{determinant} the polynomial $\det (f_i^{a_j})$ is a nonzero homogeneous 
element of $U$, which clearly has coefficients in the prime subfield of $k$.  Thus 
$p_1^{-1}, p_2^{-1}, \dots p_t^{-1}$ satisfy some nonzero non-homogeneous relation 
with coefficients in the prime subfield of $k$, contradicting the hypothesis on the 
$\{p_i\}$.

Thus the set $\mc{C}$ must be in general position, and so critically dense.  Then 
certainly $c$ must also have infinite order under $\varphi$, so that 
Hypothesis~\ref{distinct} holds.  Now $R(\varphi,c)$ is noetherian by Theorem~\ref{R 
noetherian}.
\end{proof}

Next, let us show that for generic choices of $\varphi$ and $c$ (in the sense of 
Definition~\ref{generic}), the ring $R(\varphi,c)$ is noetherian.  Because of 
Lemma~\ref{switch sides}(2), for every fixed $c \in \mb{P}^t$ we get the same class of rings $\{ 
R(\varphi,c) \mid \varphi \in \aut \mb{P}^t \}$.  Thus we might as well fix some 
arbitrary $c$ and vary $\varphi$ only.  
\begin{theorem}
\label{crit dense generic} Assume that the base field $k$ is uncountable.  Fix $c \in 
\mb{P}^t$.  There is a generic subset $Y$ of $X = \aut \mb{P}^t$ such that 
$R(\varphi, c)$ is noetherian for all $\varphi \in Y$.
\end{theorem}
\begin{proof}
By Lemma~\ref{switch sides}(2) there is no harm in assuming that $c = (1:1: \dots 
:1)$. Choose some homogeneous coordinates $(z_{ij})_{0 \leq i,j \leq t}$ for $X 
\subseteq \mb{P}(M_{t+1}(k))$.  Just as in the proof of Theorem~\ref{diag case}, we 
see that $\mc{C} = \{c_i\}_{i \in \mb{Z}}$ fails to be in general position if and 
only if there exists some $d \geq 1$ and some choice of $N = \binom{t+d}{d}$ 
nonnegative integers $a_1 < a_2 < \dots < a_N$ such that $\det f_i(c_{a_j}) = 0$, 
where the $f_i$ are the degree $d$ monomials in $U$.  

Each condition $\det f_i(c_{a_j})= 0$ is a closed condition in the coordinates of 
$X$; moreover it does not hold identically, otherwise for no choice of $\varphi$ 
would $\mc{C}$ be in general position, in contradiction to the proof of Theorem~\ref{diag case}.  
There are countably many such conditions, and so the complement $Y$ of the union of 
all of these closed subsets is generic by definition.  Thus for $\varphi \in Y$ we 
have that $\mc{C}$ is in general position and so $R(\varphi,c)$ is noetherian, by 
Theorem~\ref{R noetherian}.
\end{proof}

It is not hard, in contrast to the preceding theorems, to come by examples of 
$(\varphi,c)$ for which the $c_i$ are distinct but not even dense, much less 
critically dense.  One such example should suffice to illustrate this situation.  
\begin{example}
\label{block of size 3 not dense} Suppose that $t = 2$ and $\cha k = 0$.  Let  
\[
\varphi = \begin{bmatrix}
  1 & 1 & 0 \\
  0 & 1 & 1 \\
  0 & 0 & 1 
\end{bmatrix}
\] 
with $c = (0:0:1)$.  Then $\mc{C}$ is not dense in $\mb{P}^2$. \\
\end{example}
\begin{proof}
One easily calculates the formula $c_{-n} = (n(n-1)/2: n :1)$ for $n \in \mb{Z}$.  
Since $\cha k = 0$, the $c_i$ are obviously distinct.  But the polynomial $f = x_0x_2 
+ \frac{1}{2}x_2x_1 - \frac{1}{2}x_1^2$ vanishes at $(n(n-1)/2:n:1)$ for every $n \in 
\mb{Z}$, and so $\{ c_i \}_{i \in \mb{Z}}$ is not dense.  
\end{proof}
A similar argument will show more generally that if the Jordan canonical form of a 
matrix representing $\varphi$ has a Jordan block of size $\geq 3$ or more than one 
Jordan block of size $2$, then given any $c \in \mb{P}^t$, the set $\mc{C}$ is not 
dense in $\mb{P}^t$.  See \cite{Rothesis} for further details.  

\subsection{Improvements in characteristic zero}
Theorems~\ref{diag case} and \ref{crit dense generic} hold for an algebraically 
closed field of arbitrary characteristic.  In case where $\cha k = 0$, we will show 
that one can get a better result by invoking the following theorem of Cutkosky and 
Srinivas.  
\begin{theorem} \cite[Theorem 7]{CuSr}
\label{alg group theorem} Let $G$ be a connected commutative algebraic group defined 
over an algebraically closed field $k$ of characteristic $0$.  Suppose that $g \in G$ 
is such that the cyclic subgroup $H = \langle g \rangle$ is dense in $G$.  Then any 
infinite subset of $H$ is dense in $G$. 
\end{theorem}

The theorem has the following consequence.
\begin{proposition}
\label{dense equals crit dense} Let $\cha k = 0$.  Then $\mc{C} = \{ \varphi^i(c) 
\}_{i \in \mb{Z}}$ is critically dense in $\mb{P}^t$ if and only if $\mc{C}$ is 
Zariski dense in $\mb{P}^t$.
\end{proposition}
\begin{proof}  
If $\mc{C}$ is critically dense in $\mb{P}^t$, then $\mc{C}$ is of course dense in 
$\mb{P}^t$ by definition.    

Now assume that $\mc{C}$ is dense.  Choose a matrix $L \in \gl_{t+1}(k)$ to represent 
$\varphi$ (so L is unique up to scalar multiple).  Now set $V = \sum_{i \in \mb{Z}} k 
L^i \subseteq M_{t+1}(k)$.  Let $\wt{c} \in \mb{A}^{t+1}$ be a particular choice of 
coordinates for $c$, and think of elements of $\mb{A}^{t+1}$ as column vectors.  Then 
the linear evaluation map $\wt{\psi}: V \ra \mb{A}^{t+1}$ defined by $N \mapsto 
N\wt{c}$ descends to a map $\psi: \mb{P}V \ra \mb{P}^t$ which sends $\varphi^i$ to 
$\varphi^i(c)$ for all $i \in \mb{Z}$.  The map $\wt{\psi}$ must be surjective, else 
$\mc{C}$ would lie on a proper linear subspace of $\mb{P}^t$.  Note also that since 
$L$ satisfies its characteristic polynomial, $\dim_k V \leq t+1$. This forces 
$\wt{\psi}$ to be an isomorphism, and so $\psi$ is an isomorphism of projective 
spaces.  In particular, writing $H = \{ \varphi^i \}_{i \in \mb{Z}}$, we have via the 
automorphism $\psi$ that $H$ is dense in $\mb{P}V$.

Now let $G = \mb{P}V \cap \pgl_{t+1}(k)$.  Then $G$ is an algebraic group, since it 
is the closure of the subgroup $H$ of $\pgl_{t+1}(k)$ \cite[Proposition 1.3]{Borel}.  
Since any two elements of $V$ commute, $G$ is commutative.  Note also that $G$ is an 
open subset of the projective space $\mb{P}V$, so $G$ is irreducible and in 
particular connected.  Finally, we always assume that the field $k$ is algebraically 
closed, so the hypotheses of Theorem~\ref{alg group theorem} are all satisfied.  

Now $H$ is dense in $G$, so $H$ is critically dense in $G$ by Theorem~\ref{alg group 
theorem}.  Then as a subset of $\mb{P}V$, $H$ is critically dense in $\mb{P}V$.  
Finally, applying $\psi$ again we get that $\mc{C}$ is critically dense in $\mb{P}^t$.
\end{proof}
Thus in case $\cha k = 0$, the question of the noetherian property for $R(\varphi,c)$ 
reduces to the question of the density of $\mc{C} = \{ \varphi^i(c) \}_{i \in 
\mb{Z}}$, which is easy to analyze for particular choices of $\varphi$ and $c$.  in 
particular, $\mc{C}$ will be dense if and only if $c$ is not contained in a proper 
closed set $X \subsetneq \mb{P}^t$ with $\varphi(X) = X$.  Let us note some specific 
examples.  Note that part (1) of the following example is a significant improvement 
over Theorem~\ref{diag case} if the field has zero characteristic. 
\begin{example}
\label{small blocks dense} Let $\cha k = 0$. 
\begin{enumerate} 
\item \label{small blocks dense 1}
Suppose that $\varphi = \diag(1, p_1, \dots, p_t)$, and that the multiplicative 
subgroup of $k^{\times}$ generated by $p_1, p_2, \dots p_t$ is $\cong \mb{Z}^t$. Let 
$c$ be the point $(a_0: a_1: \dots: a_t)$.  Then $R(\varphi,c)$ is noetherian if and 
only if $a_i \neq 0$ for all $0 \leq i \leq t$.  
\item \label{small blocks dense 2}
Let
\[
\varphi = \begin{bmatrix}
  1&  1&  &  &  \\
   0& 1  & &  &  \\
   &  &  p_2&  &  \\
   &  &  & \dots  &  \\
   &  &  &  & p_t
\end{bmatrix}
\]
such that the multiplicative subgroup of $k^{\times}$ generated by the $p_2, \dots 
p_t$ is $\cong \mb{Z}^{t-1}$.  Let $c = (a_0:a_1:\dots:a_t) \in \mb{P}^t$.  Then 
$R(\varphi,c)$ is noetherian if and only if $a_i \neq 0$ for all $1 \leq i \leq t$. 
\end{enumerate}
\end{example}
\begin{proof}
(1)  Let $\phi$ be the automorphism of $U$ corresponding to $\varphi$; explicitly (up 
to scalar multiple), $\phi(x_i) = p_i x_i$, if we set $p_0 = 1$.  Suppose that $J$ is 
a graded ideal of $U$ with $\phi(J) = J$.  Then if we choose $m \gg 0$ such that $J_m 
\neq 0$, then there is some $0 \neq f \in J_m$ with $\phi(f) \in kf$, since the 
action of $\phi$ on the finite dimensional vector space $J_m$ has an eigenvector.  If 
$f = \sum b_I x_I$ (where $I$ is a multi-index), then $\phi(f) = \sum b_I p_I x_I$.  
The hypothesis on the $p_i$ forces $p_I$ to be distinct for distinct multi-indices 
$I$ of degree $m$, so $f$ must be a scalar multiple of a single monomial in the 
$x_i$.  Thus any closed set $X \subsetneq \mb{P}^t$ with $\varphi(X) = X$ is 
contained in the union of hyperplanes $\cup_{i = 0}^{t} \{x_i = 0 \}$.  It follows 
that if all $a_i \neq 0$ then $\mc{C}$ is dense.  Conversely, if some $a_i = 0$ then 
$\mc{C}$ is contained in the hyperplane $x_i = 0$ and $\mc{C}$ is not dense.  
Now the result follows from Lemma~\ref{dense equals crit dense} and Theorem~\ref{R noetherian}.  

(2) The automorphism $\phi$ of $U$ corresponding to $\varphi$ is given by $\phi(x_0) 
= x_0 + x_1$, $\phi(x_1) = x_1$, and $\phi(x_i) = p_ix_i$ for $2 \leq i \leq t$.  If 
$J$ is a graded ideal of $U$ with $\phi(J) = J$, then as above there is some $0 \neq 
f \in U_m$ with $\phi(f) \in kf$.  We leave it to the reader to show that this forces 
$f$ to be scalar multiple of a monomial in $x_1, x_2, \dots x_t$ only; the rest of 
the proof is as in part (1).  
\end{proof}

Lemma~\ref{dense equals crit dense} and Example~\ref{small blocks dense}(1) fail in 
positive characteristic.  The next example, which we thank Mel Hochster for 
suggesting, shows this explicitly.    
\begin{example}
\label{wierd ex} Let $k$ have characteristic $p > 0$ and let $y \in k$ be 
transcendental over the prime subfield $\mb{F}_p$.  Suppose that $t = 2$, and let 
$\varphi = \diag (1, y, y+1)$ and $c = (1:1:1)$.  The multiplicative subgroup of 
$k^{\times}$ generated by $y$ and $y +1$ is isomorphic to $\mb{Z}^2$, but both sets 
of points $\{c_i\}_{i \geq 0}$ and $\{c_i\}_{i \leq 0}$ are dense but not critically 
dense in $\mb{P}^t$.  The ring $R(\varphi,c)$ is neither left nor right noetherian. 
\end{example}
\begin{proof} 
It is easy to see since $y$ is transcendental over $\mb{F}_p$ that the multiplicative 
subgroup of $k^{\times}$ generated by $y$ and $y +1$ is isomorphic to $\mb{Z}^2$.

We have $c_{-n} = \varphi^n(c) = (1, y^n, (y+1)^n)$, so $c$ has infinite order under 
$\varphi$.    If $n = p^j$ for some $j \geq 0$, then $(y+1)^n =y^n + 1$.  Therefore 
$c_{-n}$ is on the line $X = \{ x_0 + x_1  - x_2 = 0 \}$ for all $n = p^j$.  On the 
other hand, suppose that $n \geq 0$ is not a power of $p$. Then some binomial 
coefficient $\binom{n}{i}$ with $0 < i < n$ is not divisible by $p$, and the binomial 
expansion of $(y+1)^n$ contains the nonzero term $\binom{n}{i}y^i$.  Since $y$ is 
transcendental over $\mb{F}_p$, this implies $(y+1)^n \neq y^n + 1$ and so $c_{-n}$ 
is not on the line $X$.  Thus for $n \geq 0$, $c_{-n}$ is on $X$ if and only if $n$ 
is a power of $p$.  It follows that the set of points $\{ c_i \}_{i \leq 0}$ is not critically 
dense.  By Theorem~\ref{R noetherian}, $R$ is not right noetherian.

Put $\mc{D} = \{c_i\}_{i \leq 0}$, and consider the Zariski closure 
$\overline{\mc{D}}$ of this set of points.  Since the line $X$ contains infinitely 
many points of $\mc{D}$, $X \subseteq \overline{\mc{D}}$.  For all $n \in \mb{Z}$, 
$\varphi^n(X)$ also contains infinitely many points of $\mc{D}$, and so $\bigcup_{n 
\in \mb{Z}} \varphi^n(X) \subseteq \overline{\mc{D}}$.  Finally, one checks that the 
$\varphi^n(X)$ are distinct lines for all $n \in \mb{Z}$.  It follows that 
$\overline{\mc{D}} = \mb{P}^2$, and $\mc{D}$ is a Zariski dense set.

Analogously, one may check that if $Y$ is the curve $\{x_1x_2 + x_2x_0 - x_0 x_1 \} = 
0$, then $c_n \in Y$ for $n \geq 0$ if and only if $n = p^j$ for some $j \geq 0$, so 
that $\{c_i\}_{i \geq 0}$ is not critically dense.  Yet a similar proof as above 
shows that $\{c_i\}_{i \geq 0}$ is 
Zariski dense.  According to Theorem~\ref{R 
noetherian}, $R$ is also not left noetherian. 
\end{proof}

\section{Algebras generated by Eulerian derivatives}
\label{eul der} 
The original motivation for our study of the algebras 
$R(\varphi,c)$ came from the results of Jordan \cite{Jordan} on algebras generated by 
two Eulerian derivatives. In this final section we show that Jordan's examples are 
special cases of the algebras $R(\varphi,c)$, and so we may use our previous results 
to answer the main open question of \cite{Jordan}, namely whether algebras generated 
by two Eulerian derivatives are ever noetherian. In fact we will prove that an 
algebra generated by a generic finite set of Eulerian derivatives is noetherian. 

Fix a Laurent polynomial algebra $k[y^{\pm 1}] = k[y, y^{-1}]$ over the base field 
$k$.  
\begin{definition}
\label{Eulder-index} For $p \in k \setminus \{0,1\}$, we define the operator $D_p \in 
\ndo_k k[y^{\pm 1}]$ by the formula $f(y) \mapsto \D \frac{f(py) -f(y)}{py-y}$.  For 
$p = 1$, we define $D_1 \in \ndo_k k[y^{\pm 1}]$ by the formula $f \mapsto  df/dy$.  
For any $p \neq 0$, we call $D_p$  an \emph{Eulerian Derivative}.
\end{definition}
It is also useful to let $y^{-1}$ be notation for the operator $y^{-1}: y^i \mapsto 
y^{i-1}$ for $i \in \mb{Z}$.
  
We now consider algebras generated by a finite set of Eulerian derivatives.  There are 
naturally two cases, depending on whether $D_1$ is one of the generators.

\begin{theorem}
\label{Eulerian derivatives 1} Suppose that $\{p_1, \dots p_t \} \in k \setminus 
\{0,1\}$ are distinct, and assume that the multiplicative subgroup of $k^{\times}$ 
these scalars generate is $\cong \mb{Z}^t$.  Let $R = k \langle D_{p_1}, D_{p_2}, \dots, 
D_{p_t} \rangle$.  Then $R \cong R(\varphi, c)$ for 
\[
\varphi = \diag(1, p_1^{-1}, p_2^{-1}, \dots, p_t^{-1})\ \text{and}\  c = (1:1: \dots 
:1).
\]
$R$ is noetherian if either $\cha k = 0$ or if the $\{p_i\}$ are algebraically 
independent over the prime subfield of $k$.
\end{theorem}
\begin{proof}
Set $p_0 = 1$ and let $w_i =y^{-1} + (p_i -1) D_{p_i}$ for $0 \leq i \leq t$.  The 
automorphism $\phi$ of $U$ corresponding to $\varphi$ is given (up to scalar 
multiple) by $\phi: x_i \mapsto p_i^{-1} x_i$ for $0 \leq i \leq t$.  An easy 
calculation shows that $S(\varphi)$ has relations $\{x_jx_i - p_j^{-1} p_i x_i x_j\}$ 
for $0 \leq i < j \leq t$; clearly these relations generate the ideal of relations 
for $S(\varphi)$, since $S(\varphi)$ has the Hilbert function of a polynomial ring in 
$t+1$ variables.

As in \cite[Section 2]{Jordan}, it is straightforward to prove the identities $w_j 
w_i - p_j^{-1} p_i w_i w_j$ for all $0 \leq i, j \leq t$, so there is a surjection of 
algebras given by  
\begin{eqnarray*}
\psi: S(\varphi) & \lra & k \langle y^{-1}, D_{p_1}, \dots, D_{p_t} 
\rangle \subseteq \ndo k[y^{\pm 1}] \\
x_i & \mapsto & w_i, \ \ \ \  0 \leq i \leq t. 
\end{eqnarray*}

Now the hypothesis that the $\{p_i\}$ generate a rank $t$ subgroup of $k^{\times}$ 
ensures that $\psi$ is injective:  this is proved for the case $t=2$ in 
\cite[Proposition~1]{Jordan}; the proof in general is analogous.  Thus $\psi$ is an 
isomorphism of algebras.  Then one checks that the image under $\psi$ of the 
subalgebra $R(\varphi,c)$ of $S(\varphi)$ is $k \langle D_{p_1}, D_{p_2}, \dots, 
D_{p_t} \rangle = R$.  

The noetherian property for $R$ follows from Example~\ref{small blocks dense}(1) in 
case $\cha k = 0$, or from Theorem~\ref{diag case} if the the $\{p_i\}$ are 
algebraically independent over the prime subfield of $k$.
\end{proof}
The case where $D_1$ is one of the generators is very similar, and we only sketch the 
proof. 
\begin{theorem}
\label{Eulerian derivatives 2} Assume that $\cha k = 0$.  Let $\{p_2, p_3, \dots p_t 
\} \in k \setminus \{0,1\}$ be distinct, and assume that the multiplicative subgroup 
of $k^{\times}$ that the $\{p_i\}$ generate is $\cong \mb{Z}^{t-1}$.  Let $R =k \langle D_1, 
D_{p_2}, D_{p_3}, \dots, D_{p_t} \rangle$.  Then $R \cong R(\varphi, c)$ for
\[
\varphi = \begin{bmatrix}
  1 & 1 &  &  &  \\
  0 &  1  & &  &  \\
   &  &  {p_2}^{-1} &  &  \\
   &  &  & \dots &  \\
   &  &  &  & p_t^{-1}
\end{bmatrix} \text{and}\ c = (0:1:1: \dots: 1).  
\]
The ring $R$ is noetherian.
\end{theorem}

\begin{proof}
Let $p_1 = 1$, and let $w_i =y^{-1} + (p_i -1) D_{p_i}$ for $1 \leq i \leq t$.
As in the preceding proposition, one calculates the relations for the algebra 
$S(\varphi)$, and using these and the identities in \cite[section 2]{Jordan}, one 
gets an algebra surjection 
\begin{eqnarray*}
\psi: S(\varphi) & \lra & k \langle y^{-1}, D_1, D_{p_2}, \dots, D_{p_t} \rangle \\
x_0 & \mapsto & -D_1 \\
x_i & \mapsto & w_i, \ \ \ \  1 \leq i \leq t. 
\end{eqnarray*}
The hypothesis on the $\{p_i\}$ implies that $\psi$ is an isomorphism, by an analogous 
proof to that of \cite[Proposition~1]{Jordan}.  Then $R(\varphi,c)$ is mapped 
isomorphically onto $R$.  The noetherian property for $R$ follows from 
Example~\ref{small blocks dense}(2). 
\end{proof}

The results above easily imply that a ring generated by a generic set of Eulerian 
derivatives is noetherian.
\begin{theorem}
\label{eul der noeth} Assume that $k$ is uncountable.  Let $V_i$ be the closed set 
$\{y_i = 0 \}$ in $\mb{A}^t$.  There is a generic subset $Y \subseteq \mb{A}^t 
\setminus \cup_{i=1}^t V_i$ such that if $(p_1, p_2, \dots, p_t) \in Y$ then $R = 
k\langle D_{p_1}, D_{p_2}, \dots D_{p_t} \rangle$ is noetherian.  
\end{theorem}
\begin{proof}
Let $k[y_1, y_2, \dots y_t]$ be the coordinate ring of $\mb{A}^t$, and write $V(f)$ 
for the vanishing set in $\mb{A}^t$ of $f \in k[y_1, y_2, \dots y_t]$.  Let $\mb{F}$ 
be the prime subfield of $k$, and set $A =\mb{F}[y_1, y_2, \dots y_t]$.  The set $Y$ 
of points $(p_1, p_2, \dots, p_t) \subseteq \mb{A}^t$ where the $\{p_i\}$ are 
algebraically independent over $\mb{F}$ is the complement in $\mb{A}^t$ of 
$\bigcup_{f \in A} V(f)$.  But since $\mb{F}$ is countable, $A$ is also countable and 
so $Y$ is a generic subset of $\mb{A}^t$ (Definition~\ref{generic}).  Now apply 
Theorem~\ref{Eulerian derivatives 1}.
\end{proof}

We can also produce an example of a ring generated by Eulerian derivatives that is 
\emph{not} noetherian, the existence of which was also an open question in \cite{Jordan}.
\begin{proposition}  Assume that $\cha k = p > 0$ and that $k$ has transcendence degree at least $1$ over
its prime subfield $\mb{F}_p$.  Then there exist scalars $p_1, p_2 \in k$ such that 
the ring $k \langle D_{p_1}, D_{p_2} \rangle$ is not noetherian.
\end{proposition}
\begin{proof}
let $y \in k$ be transcendental over $\mb{F}_p$.  Consider the ring $R(\varphi,c)$ of 
Example~\ref{wierd ex}, where $\varphi = \diag(1, y, y+1)$ and $c = (1 : 1: 1)$.  As 
in Example~\ref{wierd ex}, the scalars $y, y+1$ generate a rank $2$ multiplicative 
subgroup of the field $k$, so setting $p_1 = y$ and $p_2 = y+1$ we have $R(\varphi, 
c) \cong k \langle D_{p_1}, D_{p_2} \rangle$ by Theorem~\ref{Eulerian derivatives 
1}.  But as we saw in Example~\ref{wierd ex}, $R(\varphi,c)$ is not noetherian.
\end{proof} 

%___________________________________________________________________________________
\appendix
\section{Castelnuovo-Mumford regularity}
In this appendix, we discuss the notion of Castelnuovo-Mumford regularity and use 
some recent results in this subject to prove the technical lemmas in the main body of 
the paper.  We thank Jessica Sidman for alerting us to these methods.

Let $U = k[x_0,x_1, \dots x_t]$ be a polynomial ring over an algebraically closed 
field $k$, graded as usual with $\deg(x_i) = 1$ for all $i$.  All multiplication in 
this appendix is commutative, and so we omit the $\circ$ notation which we introduced in 
\S\ref{Zhang twists}.  Generally speaking, the notion of regularity for a $U$-module 
$M$ is a convenient way of encapsulating information about the degrees of the 
generators of all of the syzygies of $M$.  
\begin{definition} \cite[page 509]{E}
\label{regular-index} \label{defreg} Let $M \in U\gr$.  Take a minimal graded free 
resolution of $M$: 
\[
0 \ra \bigoplus_{i=1}^{r_{(t+1)}} U[-e_{i,t+1}] \ra \dots \ra \bigoplus_{i=1}^{r_0} 
U[-e_{i,0}] \ra M \ra 0.
\]
If $e_{i,j} \leq m + j$ for all $i, j$ then we say that $M$ is \emph{$m$-regular}. 
The \emph{regularity} of $M$, $\reg M$, is the smallest integer $m$ for which $M$ is 
$m$-regular (if $M = 0$ then we set $\reg M = - \infty$).
\end{definition}
There are other equivalent characterizations of regularity, with different 
advantages; see for example \cite[Definition 3.2]{BaMu}.  

Regularity behaves well with respect to exact sequences.
\begin{lemma} \cite[Corollary~20.19]{E}
\label{reg ex seq} let $0 \ra M' \ra M \ra M'' \ra 0$ be a short exact sequence in 
$U\gr$.  Then 
\begin{enumerate}
\item \label{reg ex seq 1}
$\reg M' \leq \max(\reg M, \reg M'' +1)$
\item \label{reg ex seq 2} 
$\reg M \leq \max(\reg M', \reg M'')$  
\item \label{reg ex seq 3}
$\reg M'' \leq \max(\reg M' -1, \reg M)$. 
\hfill $\Box$
\end{enumerate} 
\end{lemma}

Let us define some related notions.  For $I$ a graded ideal of $U$, we define the 
\emph{saturation} of $I$ to be 
\[
I^{sat} = \{ x \in U | (U_{\geq n}) x \subseteq I\ \text{for some}\ n \}. 
\]
The ideal $I^{sat}$ is the unique largest extension of $I$ inside $U$ by a torsion 
module.  If $I^{sat} = I$ then we say that $I$ is \emph{saturated}.  

A module or ideal that is $m$-regular stabilizes in degree $m$ in the following important ways.    
\begin{lemma}
\label{reg easy facts} 
\begin{enumerate}
\item \label{reg easy facts 1}
If $M \in U\gr$ is $m$-regular, then $M$ is generated in degrees less than or equal to $m$. 
\item \label{reg easy facts 2} 
If $I$ is a graded ideal of $U$ then $(\sat I)_{\geq m} = I_{\geq m}$ for $m = \reg I$. 
\item \label{reg easy facts 3} 
If $I$ is a graded ideal of $U$, then $I$ is $m$-regular if and only if $I_{\geq m}$ 
is $m$-regular. 
\item \label{reg easy facts 4} 
If $M \in U\gr$ then the Hilbert function of $F(n) = \dim_k M_n$ of 
$M$ agrees with the Hilbert polynomial of $M$ in degrees $\geq (\reg M +1) - (t+1 - 
\pd(M))$, where $\pd(M)$ is the projective dimension of $M$.
\end{enumerate}
\end{lemma}
\begin{proof}
(1) is immediate from Definition~\ref{defreg}, and (2) and (3) follow from \cite[Definition 3.2]{BaMu}.
For (4), note that the Hilbert function of $U$, $f(n) = \dim_k U_n$, agrees with its Hilbert 
polynomial $(n+t)(n+t-1) \dots (n+1)/{t !}$ for $n \geq -t$.  Then calculating the 
Hilbert function of $M$ from its minimal free resolution, the result follows from the 
definition of regularity.  
\end{proof}

The regularity of an ideal $I \subseteq U$ might be much greater than the minimal 
generating degree of $I$, but at least there is the following bound.   
\begin{lemma} \cite[Proposition~3.8]{BaMu}
\label{exp bound} Let $I$ be a homogeneous ideal of $U$, and let $d$ be the minimal generating degree
of $I$.  Then $\reg I \leq (2d)^{t !}$. \hfill $\Box$
\end{lemma} 

The key ingredients in the proofs of our needed lemmas will be
the following recent theorems of Conca and Herzog concerning the regularity of products.

\begin{theorem} 
\label{CoHe 1} \cite[Theorem~2.5]{CoHe} If I is a graded ideal of $U$ with $\dim U/I 
\leq 1$, then for any $M \in U\gr$ we have $\reg IM \leq \reg I + \reg M$.  \hfill 
$\Box$
\end{theorem}

\begin{theorem}
\label{CoHe 2} \cite[Theorem~3.1]{CoHe} Let $I_1$, $I_2$, \dots, $I_e$ be (not 
necessarily distinct) nonzero ideals of $U$ generated by linear forms.  Then $\reg 
(I_1 I_2 \dots I_e) = e$.  \hfill $\Box$
\end{theorem}

We may now make the following observations about the regularity of an ideal of a finite set of points with 
multiplicites.
\begin{lemma} 
\label{reg results}
\label{better reg results}
Let $d_1, d_2, \dots, d_n$ be distinct points of $\mb{P}^t$ with ideals $\m_1, \m_2, \dots, \m_n$.
Let $0 < e_i$ for $1 \leq i \leq n$ and set $e = \sum e_i$.  Let $J = \m_1^{e_1} \cap \m_2^{e_2} \cap \dots \cap \m_n^{e_n}$.
\begin{enumerate}
\item 
$\reg J \leq e$. 
\item  
If the points $d_1, d_2, \dots, d_n$ do not all lie on a line, then $\reg J \leq e-1$.  
\end{enumerate} 
\end{lemma}
\begin{proof}
(1) Let $I =\m_1^{e_1} \m_2^{e_2} \dots \m_n^{e_n}$. 
It is easy to see that $J$ is saturated, and that $J/I$ is torsion, so that $J$ is the saturation of $I$.  
By  Theorem~\ref{CoHe 2}, $\reg I \leq e$.  Then $J_{\geq e} = (I^{sat})_{\geq e} = I_{\geq e}$ by
Lemma~\ref{reg easy facts}(2).  By Lemma~\ref{reg easy facts}(3), $\reg J \leq e$.

(2) The hypothesis on the points forces some three of the points $\{d_i\}$ to be 
non-collinear (in particular $n \geq 3$); by relabeling we may assume that $d_1, d_2, 
d_3$ do not lie on a line.  Then one may check that $U/(\m_1 \cap \m_2 \cap \m_3)_{\geq 1}$ 
is isomorphic to $(U/\m_1 \oplus U/\m_2 \oplus U/\m_3)_{\geq 1}$ and so this module is $1$-regular.
Then by Lemma~\ref{reg easy facts}(2), $U/(\m_1 \cap \m_2 \cap \m_3)$ is $1$-regular, so
using Lemma~\ref{reg ex seq} we get that $\m_1 \cap \m_2 \cap \m_3$ is $2$-regular.

Now $L = (\m_1 \cap \m_2 \cap \m_3)(\m_1^{e_1 -1} \cap \m_2^{e_2 -1} \cap \m_3^{e_3 
-1} \cap \m_4^{e_4} \cap \dots \cap \m_n^{e_n})$ is $(e-1)$-regular, using 
part (1) and Theorem~\ref{CoHe 1}.  Since $L^{sat} = J$, a similar argument as in part (1) 
shows that $J$ is $(e-1)$-regular. 
\end{proof}

Now we prove the series of results which we used in the body of the paper.  The first is an immediate corollary of
the proof of part (1) of the preceding lemma.
\begin{lemma} (Lemma \ref{products versus intersections})
\label{products versus intersections app} Let $m_1, m_2, \dots, m_n$ be the ideals of 
$U$ corresponding to distinct points $d_1, \dots d_n$ in $\mb{P}^t$.  Then 
$(\prod_{i=1}^n \m_i)_{\geq n} = (\bigcap_{i=1}^n \m_i)_{\geq n}$. \hfill $\Box$
\end{lemma}

Next is a simple Hilbert function calculation, which is presumably well known.  Because 
this lemma is so fundamental above, we include a brief proof which uses the methods 
of regularity.
\begin{lemma} (Lemma \ref{Hilbert of point ideal})
\label{Hilbert of point ideal app} Let $e_i > 0$ for all $1 \leq i \leq n$, and let 
$e = \sum e_i$.  Set $J = \bigcap_{i=1}^n \m_i^{e_i}$ for some distinct point ideals 
$\m_1, \m_2, \dots, \m_n$. Then $\dim_k J_m = 
\binom{m+t}{t} - \sum_i \binom{e_i+t-1}{t}$ for all $m \geq e -1$. \\
In particular, if $J = \bigcap_{i=1}^n \m_i$ then $\dim_k J_m = \binom{m+t}{t} - n$ 
for $m \geq n-1$. 
\end{lemma}
\begin{proof}
We leave it to the reader to show that the Hilbert polynomial of the module $U/\m_i^{e_i}$ 
is the constant $\binom{e_i+t-1}{t}$, for example by induction on $e_i$.  Then since the points $d_i$ are distinct, 
the Hilbert polynomial of $J$ is $H(m) = \binom{m+t}{t} - \sum_i 
\binom{e_i+t-1}{t}$. 

By Lemma~\ref{reg results}(1) $\reg J \leq e$.  By the Auslander-Buchsbaum formula we 
have since $\dep(U/J) =1$ that $\pd(U/J) = t$, and so $\pd(J) = t-1$.  Now it follows 
from Lemma~\ref{reg easy facts}(4) that the Hilbert function of $J$ agrees with its 
Hilbert polynomial in degrees $\geq e-1$. 
\end{proof}

\begin{lemma} (Lemma~\ref{case not on a line})
\label{case not on a line app} Let the points $d_1, d_2, \dots, d_n, d_{n+1} \in 
\mb{P}^t$ be distinct, and assume that the points $d_1, \dots, d_n$ do not all lie on 
a line.  Let $\m_i \subseteq U$ be the homogeneous ideal corresponding to $d_i$.  
\begin{enumerate}
\item \label{case not on a line app 1}
$(\bigcap_{i=1}^n \m_i)_{n-1} (\m_{n+1})_1 = (\bigcap_{i = 1}^{n+1} \m_i)_n$. 
\item \label{case not on a line app 2}
$(\bigcap_{i=1}^n \m_i)_{n-1} (\m_1)_1 = (\bigcap_{i = 2}^n \m_i \cap \m_1^2)_n$. 
\item \label{case not on a line app 3}
$(\bigcap_{i=2}^n \m_i \cap \m_1^2)_n (\m_{n+1})_1 = (\bigcap_{i = 2}^{n+1} \m_i \cap 
\m_1^2)_{n+1}$. 
\item \label{case not on a line app 4}
Let $b_1, b_2 \in \mb{P}^t$, with corresponding ideals $\mf{n_1}, \mf{n_2}$, be such 
that $b_j \neq d_i$ for $j = 1, 2$ and $1 \leq i \leq n$.  Then $(\bigcap_{i=1}^n 
\m_i \cap \mf{n}_1)_n = (\bigcap_{i=1}^n \m_i \cap \mf{n}_2)_n$ implies $b_1 =b_2$.  
\end{enumerate} 
\end{lemma} 
\begin{proof}
(1) Set $K = \bigcap_{i=1}^n \m_i$, $L = \m_{n+1}$, and $M = \bigcap_{i = 1}^{n+1} 
\m_i$.  By Lemma~\ref{reg results}, we have that $\reg 
K \leq n-1$ and $\reg L \leq 1$.  Then by Theorem~\ref{CoHe 1}, $\reg (KL) \leq n$.  
Since clearly $M = (KL)^{sat}$, by Lemma~\ref{reg easy facts}(2) it follows that 
$(KL)_n = M_n$.  Finally, by Lemma~\ref{reg easy facts}(1), $K$ is generated in 
degrees $\leq n-1$ and $L$ is generated in degree $1$.  Thus $K_{n-1} L_1 = (KL)_n = 
M_n$.  

(2)-(3) The proofs of these parts are very similar to the proof of (1) and are 
omitted.  

(4) The ideals $K = ( \bigcap_{i=1}^n \m_i \cap \mf{n}_1)$ and $L = (\bigcap_{i=1}^n 
\m_i  \cap \mf{n}_2)$ are each $n$-regular by Lemma~\ref{better reg results}, so both 
are generated in degrees $\leq n$.  Now since $b_1 \neq d_i$ for all $i$, if $b_1 
\neq b_2$ then the ideals $K$ and $L$ must differ in large degree, so they must 
differ in degree $n$.
\end{proof}

We close with a simple application of regularity to the analysis of bounds for 
$\uext$ groups.  

\begin{lemma}
\label{reg ext app} Given $M, N \in U\gr$, there is a constant $d \in \mb{Z}$, 
depending only on $\reg M$ and $\reg N$, such that $\reg (\uext^i_U(M,N)) < d$ for 
all $i \geq 0$. 
\end{lemma}
\begin{proof}
Take a minimal graded free resolution of $M$:
\[
0 \ra \bigoplus_{i=1}^{r_{(t+1)}} U[-e_{i,t+1}] \ra \dots \ra \bigoplus_{i=1}^{r_0} 
U[-e_{i,0}] \ra M \ra 0.
\]
By the definition of regularity, $e_{i,j} \leq \reg M + (t+1)$ for all $i, j \geq 
0$.  Applying $\uhom(- ,N)$ to the complex with the $M$ term deleted produces a 
complex 
\[
0 \ra L_0 \overset{\psi^0}{\ra} L_1 \overset{\psi^1}{\ra} \dots \overset{\psi^t}{\ra} 
L_{t+1} \ra 0 
\]
where $L_j = \bigoplus_{i=1}^{r_j} N[-e_{i,j}]$. Then $\reg L_j \leq (\reg M + \reg N 
+ t +1)$ for all $j \geq 0$.  

Now consider the map $\psi^i: L_i \ra L_{i+1}$ for some $i \geq 0$.  Certainly $L_i$ 
is generated in degrees less than or equal to $\reg L_i$, by Lemma~\ref{reg easy 
facts}(1).  Then $\im \psi^i$ is also generated in degrees less than or equal to 
$\reg L_i$.  By Lemma~\ref{exp bound}, $\reg (\im \psi^i) \leq f(\reg L_i)$ where 
$f(x) = (2x)^{t !}$.  By Lemma~\ref{reg ex seq}(1), 
\[
\reg (\ker \psi^i) \leq \max (\reg L_i, \reg (\im \psi^i) +1) \leq f(\reg L_i) + 1. 
\]
Finally, $\uext^i(M,N) \cong \ker \psi^i/\im \psi^{i-1}$ and so by \ref{reg ex 
seq}(3),
\begin{gather*}
\reg (\uext^i(M,N)) \leq \max(f(\reg L_i)+1, f(\reg L_{i-1})) \\
\leq f(\reg M + \reg N + t + 1) +1
\end{gather*}
and thus we may take $d = f(\reg M + \reg N + t + 1) +1$.
\end{proof}

The final lemma is an easy corollary of the preceding one.
\begin{lemma} (Lemma~\ref{comm ext lemma})
\label{comm ext lemma app} Let $I,J$ be any homogeneous ideals of $U$, and let $\phi$ 
be an automorphism of $U$. Then there is some fixed $d \geq 0$ such that for all $n 
\in \mb{Z}$ such that $U/(I + \phi^n(J))$ is bounded, $\uext_U^i(U/I, U/\phi^n(J))$ 
is also bounded for all $i$, with $d$ as a right bound.  
\end{lemma}
\begin{proof}
If $U/(I + \phi^n(J))$ is bounded, then $E^n = \uext_U^i(U/I, U/\phi^n(J))$ is 
certainly also bounded, since it is killed by $I + \phi^n(J)$.  It is clear from the 
definition of regularity that the modules $\{ U/\phi^n(J) \}_{n \in \mb{Z}}$ all have 
the same regularity, and so by Lemma~\ref{reg ext app} there is some bound $d \geq 0$ 
such that $\reg E^n \leq d$ for all $n \in \mb{Z}$.  Then if $E^n$ is bounded, $d$
 is a right bound for it,  by Lemma~\ref{reg easy facts}(4).  
\end{proof}

\providecommand{\bysame}{\leavevmode\hbox to3em{\hrulefill}\thinspace}

\end{document}